\input amstex
\documentstyle {amsppt} 
\magnification=1200
\vsize=9.5truein
\hsize=6.5truein
\nopagenumbers 
\nologo


\topmatter

\title
Representation zeta functions
of wreath products with finite groups
\endtitle

\leftheadtext{Representation zeta functions
of wreath products}
\rightheadtext{Representation zeta functions
of wreath products}

\author
Laurent Bartholdi and Pierre de la Harpe
\endauthor

\address
\hskip-.6cm
Laurent Bartholdi,  
Mathematisches Institut,
Bunsenstra\ss e 3--5,
Georg August-Universit\"at,
D-37073 G\"ottingen.
Email:  laurent.bartholdi\@gmail.com
\endaddress

\address
\hskip-.6cm
Pierre de la Harpe, Section de Math\'ematiques, 
Universit\'e de Gen\`eve,
C.P. 64, 
\newline
CH--1211 Gen\`eve 4.
Email: Pierre.delaHarpe\@unige.ch
\endaddress


\thanks
We are grateful to 
Slava Grigorchuk, who has encouraged one of us
to begin the work presented here, 
to Don Zagier for helpful discussions, encouragement, and calculations,
and  to the Swiss National Science Foundation,
for its support to visiting colleagues.
\endthanks

\abstract
Let $G$ be a group which has a finite number $r_n(G)$ of
irreducible linear representations in $GL_n(\bold C)$ for all $n \ge 1$.
Let $\zeta(G,s) = \sum_{n=1}^{\infty} r_n(G) n^{-s}$
be its representation zeta function.
\par
First, in case $G = H \wr_X Q$ is a permutational wreath product
with respect to a permutation group $Q$ on a finite set $X$,
we establish a formula for $\zeta(G,s)$ in terms of
the zeta functions of $H$ and of subgroups of $Q$, 
and of the M\"obius function associated to the lattice $\Pi_Q(X)$
of partitions of $X$ in orbits under subgroups of $Q$.
\par
 
Then, we consider groups 
$W(Q,k) = ( \cdots (Q \wr_X Q) \wr_X  Q \cdots) \wr_X Q$
which are iterated wreath products (with $k$ factors $Q$), 
and several related infinite groups $W(Q)$,
including the profinite group $\varprojlim_k W(Q,k)$, 
a locally finite group $\lim_k W(Q,k)$,
and several  finitely generated dense subgroups of $\varprojlim_k W(Q,k)$.
Under convenient hypotheses (in particular $Q$ should be perfect),
we show that  $r_n(W(Q)) < \infty$ for all $n \ge 1$,
and we establish that the Dirichlet series $\zeta(W(Q),s)$ 
has a finite and positive abscissa of convergence 
$\sigma_0 = \sigma_0(W(Q))$.
Moreover, the function $\zeta(W(Q),s)$ 
satisfies a remarkable functional equation
involving $\zeta(W(Q),es)$ for $e \in \{1,\hdots,d\}$,
where $d = \vert X \vert$. 
As a consequence of this,
we exhibit some properties of the function,
in particular that $\zeta(W(Q),s)$ 
has a root--type singularity at $\sigma_0$,
with a finite value at $\sigma_0$ 
and a Puiseux expansion around~$\sigma_0$.
\par
We finally report some numerical computations 
for $Q = A_5$ and $Q = PGL_3(\bold F_2)$.
\endabstract



\endtopmatter

\document

\head{\bf
1.~Introduction
}\endhead

Let $G$ be a group. 
We denote by $\widehat G$ 
the set of equivalence classes of irreducible linear representations
of $G$ in complex vector spaces of finite dimension;
in case $G$ is a topological group, 
we assume that the representations are continuous.
The set $\widehat G$ is called here (and somewhat abusively)
the {\it dual} of $G$.
For $n \ge 1$, let $r_n(G) \in \{0,1,2,\hdots,\infty\}$
denote the number of $\pi \in \widehat G$ of degree $\deg \pi = n$.
The group $G$ is {\it rigid} if $r_n(G) < \infty$ for all $n \ge 1$.
\par

For example,  an infinite cyclic group is not rigid,
since $r_1(\bold Z) = \infty$;
more generally, if a finitely generated group $G$ 
has a finite index normal subgroup with infinite abelianisation, 
then $G$ is not rigid.
On the other hand, rigid groups include 
finite groups,
compact semisimple Lie groups (as a consequence of Weyl's formula 
for the dimensions of their irreducible representations),
compact groups $H(\bold Z_p)$ where
$H$ denotes a semisimple group defined over 
some field $\bold Q_p$ of $p$--adic numbers 
(Proposition 4.1 in \cite{SaAR--94}),
and many arithmetic groups \cite{LaLu--08}.
\par

The {\it representation zeta function} of a rigid group $G$ 
is the Dirichlet series
$$
\zeta_G(s) \, \Doteq \, \sum_{\pi \in \widehat G} (\deg \pi)^{-s}  
\, = \, \sum_{n \ge 1} r_n(G) n^{-s} \hskip.1cm .
\tag{1.1}
$$
For a rigid group $G$ such that $\sum_{n \ge 1} r_n(G) = \infty$,
namely such that $\zeta_G(s)$ diverges at $s=0$,
it is a classical result  (see e.g.\ Theorems 7 and 10 in \cite{HaRi--15})
that the {\it abscissa of convergence} of $\zeta_G(s)$ is given by
$$ 
\sigma_0(G) \, \Doteq \, \limsup_{n \to \infty} 
\frac{ \ln \left( \sum_{j=1}^n r_j(G) \right)}{\ln n} ,
\tag{1.2}
$$
so that it measures in some sense the growth of the $r_j(G)$'s,
and that $\sigma_0(G)$ is a singular point of the 
{\it function} $\zeta_G(s)$.
We have {\it a priori} $0 \le \sigma_0(G) \le \infty$,
and we will show that $0 < \sigma_0(G) < \infty$
for the groups which appear below.
Whenever convenient, 
we write $\zeta(G,s)$ instead of $\zeta_G(s)$.
\par
For example, $\zeta(SU(2),s)$ is the Riemann zeta function
and $\sigma_0(SU(2)) = 1$ since we have
$r_n(SU(2)) = 1$ for any $n \ge 1$.
If $G$ is a connected compact Lie group which is semisimple
and simply connected, the value of $\sigma_0(G)$ 
and information about $\zeta_G(s)$ are implicitely contained
in formulas due to Weyl (see Theorem~20).
If $G$ is the Gupta--Sidki group corresponding to an odd prime $p$,
it is a particular case of Theorem 1.3 of \cite{PaTe--96}
that $\sigma_0(G) \ge p-2$, 
but we do not know more about $\zeta_G(s)$.
We would also like to know more about $\zeta_G(s)$ and $\sigma_0(G)$
for $G$ the Grigorchuk group.
\par

On the one hand, the function $\zeta_G(s)$ captures a very small part only of the
representation theory of $G$;
when $G$ is finite, it is just a way
to organise the information contained in the {\it degree pattern} of $G$,
which is the list of the integers $\deg \pi$, including multiplicities 
(see  \cite{Hupp--98}).
On the other hand, $\zeta_G(s)$
happens to be strongly related to several interesting questions,
as shown by numerous articles including 
\cite{Zagi--96}, \cite{Lars--04}, \cite{LuMa--04}, 
\cite{LiSh--05}, \cite{Jaik--06}, \cite{LaLu--08}, and \cite{Avni}.
As far as we know, the first appearance of these zeta functions is in
\cite{Witt--91}, where his Formula (4.72) relates 
the evaluation at $2g - 2$ of $\zeta_G(s)$ 
to the volume of the moduli space 
of flat connections of $G$--principal bundles over $\Sigma_g$, 
where $G$ is a compact, simple, simply connected Lie group 
and $\Sigma_g$ an orientable closed surface of genus~${g \ge 2}$.
In general, the study of $\zeta_G$ is related to 
the subject of  \lq\lq subgroup growth\rq\rq \ \cite{LubSe--03},
which is essentially that of counting permutation representations in symmetric groups
(instead of linear representations of some sort, as here).

\bigskip

The aim of our  work  is to exhibit  properties
of representation zeta functions for groups which are wreath products,
and in particular for groups $G$ which are isomorphic 
to some wreath products with themselves, 
namely of the form $G \cong G \wr_X Q$.
Our main result is to show that there exists a group $G$
for which $\zeta_G(s)$ has a {\it root singularity}
at its abscissa of convergence (see Theorem~3),
in sharp contrast to what happens 
for compact Lie groups (see Theorem~20).
This function $\zeta_G(s)$ satisfies a remarkable {\it functional equation},
from which most of what we know about the function follows
(see Theorem~3 again).
\par

Our results are stated for profinite groups in Sections 2 and 3,
and proven in Sections 4 to~6.
Section 7 is a description of other classes of groups 
to which these results apply (see Theorem~24).
Section 8 comments on our numerical computations.
Some elementary facts on representation zeta functions
are collected in Section~9.
The final Section~10 contains a remark
concerning the \lq\lq unitary variation\rq\rq \ 
of representation zeta functions
and a question about a possible strengthening of Kazhdan's Property~(T).

\head{\bf
2.~Wreath products with a finite permutation group
}\endhead

Let $H$ be a topological group 
and let $Q$ be a finite group acting on a finite set $X$.
Let
$$
G \, = \, H \wr_X Q \, = \, B \rtimes Q ,
\hskip.5cm \text{with} \hskip.2cm B \Doteq H^X ,
\tag{2.1}
$$
denote the corresponding {\it wreath product}.
Here, $H^X$ denote the group of all applications from $X$ to $H$,
with pointwise multiplication, and $\rtimes$ indicates
the semi--direct product with respect to 
the natural action of $Q$ on $B$,
given by $(qb)(x) = b(q^{-1}x)$ 
for $q \in Q$, $b \in B$, and $x \in X$.
The topology on $B \rtimes Q$ is that for which 
$B = H^X$ has the product topology and is an open subgroup.
\par

Our first goal is to revisit part of the 
{\it representation theory of wreath products.}
This is a classical subject:
see among others \cite{Spec--33}, \cite{Kerb--71}
(which contains historical comments 
on wreath products and their representations), 
and \cite{Hupp--98, \S~25};
let us also mention that
the theory of iterated wreath products (see Section~3)
goes back to Kaloujnine 
\cite{Kalo--45, Kalo--48} and his students.
As a consequence, an important part of what follows 
consists of variations on standard themes.
This holds in particular for our first result,
before the statement of which 
we need to define the appropriate M\"obius function. 
\par

Any subgroup $S$ of $Q$ gives rise to a partition of $X$
in $S$--orbits; we denote by $\Pi_Q(X)$ 
the {\it lattice} of all partitions of $X$ of this kind,
where, for $P,P' \in \Pi_Q(X)$,  we have $P \le P'$
if $P$  is a refinement of $P'$ 
(namely if every {\it block} of $P$ is contained in a block of $P'$).
For a partition $P = (P_1,\hdots,P_\ell)$ in $\Pi_Q(X)$,
its {\it stability subgroup} is the corresponding subgroup
$$
Q_P \, = \, 
\left\{ q \in Q \hskip.1cm \vert \hskip.1cm
q(P_1) = P_1, \hdots, q(P_\ell) = P_\ell \right\} 
\tag{2.2}
$$
of $Q$.
The smallest partition $\widehat 0$, 
that for which all blocks are singletons, 
has stability subgroup the subgroup of $Q$
of trivially acting elements 
(the one--element subgroup $\{1\}$ if the action is faithful),
and the largest partition $\widehat 1$,
that for which the blocks are the $Q$--orbits,
has stability subgroup $Q$
(so that $\widehat 1$ has just one block 
if $Q$ acts transitively on~$X$).
We denote by $\mu_X$ 
the M\"obius function of the lattice $\Pi_Q(X)$;
recall that the domain of~$\mu_X$ is $\Pi_Q(X)^2$,
that $\mu_X(P,P'') = 0$ unless  $P \le P''$,
that $\mu_X(P,P) = 1$ for all $P \in \Pi_X(Q)$,
and that
$$
\mu(P,P'') \, = \, 
- \sum_{P' \in \Pi_Q(X) \atop P \le P' < P''}
\mu(P,P') 
\hskip.5cm \text{if} \hskip.2cm P < P'' .
\tag{2.3}
$$
We refer to \cite{Stan--97}; see in particular his Example 3.10.4
for the lattice of all partitions of~$X$, which is our $\Pi_Q(X)$
for $Q$ the group of all permutations of $X$.
For an easier example, see below, in the proof of Proposition~6,
Formula (4.6).
\par

Our first result is the following consequence of Clifford's theory, 
proven in Section~4.

\proclaim{1.~Theorem}
Let $H$ be a topological group, 
let $Q$ be a finite group acting on a finite set $X$,
and let $H \wr_X Q$ be the corresponding permutational wreath product, 
as above.
Then $H \wr_X Q$ is rigid if and only if $H$ is rigid.
Moreover, when this is the case, we have
$$
\aligned
\zeta(H \wr_X Q,s) \, &= \,
\sum_{P \in \Pi_Q(X)} [Q : Q_P ]^{-1-s}
\zeta_{Q_P}(s)
\\
& \hskip2cm
\sum_{P' = (P'_1,\hdots,P'_\ell) \ge P} \hskip.1cm \mu_X( P,P')  \hskip.1cm 
\zeta_H(\vert P'_1 \vert s) \cdots \zeta_H(\vert P'_\ell \vert s) 
\endaligned
\tag{2.4}
$$
(see also Formula (8.1) below).
\endproclaim

For example, when $Q$ is the permutation group of three objects
acting on $\{1,2,3\}$, we have
$$
\aligned
\zeta(H \wr_3 S_3, s) \, &= \, 
6^{-1-s} \Big( 
\zeta_H(s)^3 - 3 \zeta_H(s)\zeta_H(2s) + 2 \zeta_H(3s) \Big)
\\
\, &+ \, 
3 \times 3^{-1-s} \times 2 \Big( 
\zeta_H(s)\zeta_H(2s) -\zeta_H(3s) \Big)
\\
\, &+ \, \Big(
2 + 2^{-s} \Big) \zeta_H(3s)
\\
\, &= \, 6^{-1-s} \zeta_H(s)^3 +
\Big( - \hskip.1cm 6^{-s}/2 + 2 \times 3^{-s} \Big) 
\zeta_H(2s) \zeta_H(s)
\\
\, &+ \, 
\Big(
6^{-s}/3 - 2 \times 3^{-s} + 2^{-s} + 2 \Big) 
\zeta_H(3s) .
\endaligned
\tag{2.5}
$$
In particular,
for the Weyl group $C_2 \wr_3 S_3$ of type $B_3$ and of order $48$,
Formula (2.5) specialises to
$$
\zeta(C_2 \wr_3 S_3,s) \, = \, 
4 + 2 \times 2^{-s} + 4 \times 3^{-s} ,
\tag{2.6}
$$
so that we recover the well--known degree pattern
$1,1,1,1,2,2,3,3,3,3$.
For other specialisations, see  Sections~4 and~8.

\head{\bf
3.~Iterated wreath products 
}\endhead

A finite group $Q$ acting on a finite set $X$ gives rise 
to a tower of iterated wreath products
$$
W(Q,0) \, = \, \{1\}, \hskip.3cm
W(Q,1) \, = \, Q, \hskip.3cm
\hdots , \hskip.3cm
W(Q,k+1) \, = \, W(Q,k) \wr_X Q, \hskip.3cm
\hdots ;
\tag{3.1}
$$
observe that
$$
\vert W(Q,k) \vert \, = \, \vert Q \vert ^{(d^{k}-1)/(d-1)} .
\tag{3.2}
$$
Recall the associativity of wreath products:
given two finite groups $Q,R$
acting respectively on two finite sets $X,Y$, 
there is a natural action of $Q \wr_Y R$ on $X \times Y$
and a natural isomorphism
$$
\big(H \wr_X Q \big) \wr_Y R \, = \, 
H \wr_{(X \times Y)} \big(Q \wr_Y R) .
\tag{3.3}
$$
In particular, we also have wreath product decompositions
$$
W(Q,k+1) \, = \,  Q \wr_{X^k} W(Q,k)
\tag{3.4}
$$
and corresponding natural split epimorphisms
$$
W(Q,k+1) \, \longrightarrow \, W(Q,k) ,
\tag{3.5}
$$
which give rise to a  profinite group
$$
W^{\text{{\it prof}}}(Q) \, = \, \varprojlim_k W(Q,k) .
\tag{3.6}
$$
Note the isomorphism
$$
W^{\text{{\it prof}}}(Q) \, \approx \, W^{\text{{\it prof}}}(Q) \wr_X Q
\tag{3.7}
$$
which  plays a fundamental role below.
A rich class of groups of this kind appear in \cite{Sega--01}.
\par

For example, if $Q = S_d$ is the symmetric group
acting in a standard way on $\{1,\hdots,d\}$, 
then $W^{\text{{\it prof}}}(S_d)$ is the full automorphism group
of the infinite $d$--ary regular rooted tree
(more on these trees below in Section~7).
Groups defined by (3.6) are far from being rigid in general:
for example, if $C_p$ is the cyclic group of prime order $p$
acting on itself by multiplication, then
$$
r_n(W^{\text{{\it prof}}}(C_p)) \, = \, \left\{ \aligned
\infty \hskip.2cm &\text{if $n=p^e$ for some $e \ge 0$,}
\\
0 \hskip.2cm &\text{otherwise.}
\endaligned \right.
\tag{3.8}
$$
It follows that $W^{\text{{\it prof}}}(S_d)$ is not rigid,
since $r_n(W^{\text{{\it prof}}}(S_d)) \ge r_n(W^{\text{{\it prof}}}(C_2))$
for all $d \ge 2$ and $n \ge 1$.
\par

The situation is radically different if $Q$ is perfect.
More precisely, and this is {\it the  main purpose of this paper,} 
our goal is to show that the representation Dirichlet series 
of the group $W^{\text{{\it prof}}}(Q)$ 
has a finite abscissa of convergence
$\sigma_0 = \sigma_0( W^{ \text{{\it prof}} }(Q) ) > 0$,
and that the resulting function,  
holomorphic in the half--plane 
$\{\operatorname{Re}(s) > \sigma_0\}$, 
has  remarkable properties.

\proclaim{2.~Theorem}
Let $Q$ be a finite group acting transitively 
on a finite set $X$ of size $d \ge 2$. 
The following properties are equivalent:
\roster
\item"(i)" the finite group $Q$ is perfect,
\item"(ii)" the profinite group $W^{\text{{\it prof}}}(Q)$ is rigid.
\endroster
\endproclaim

\proclaim{3.~Theorem} Let $Q$ be as in the previous theorem
and assume that Properties (i) and (ii) hold. Then
the abscissa of convergence $\sigma_0$ 
of $\zeta(s) = \zeta(W^{\text{{\it prof}}}(Q),s)$  is a finite positive number, 
and the function $\zeta(s)$ has a singularity 
at $\sigma_0$ 
with a Puiseux expansion of the form
$$
\zeta(s) \, = \, 
\sum_{n=0}^{\infty} a_n (s-\sigma_0)^{\frac{n}{e}} 
\tag{3.9}
$$
for some integer $e$ with $2 \le e \le d$.
The function $\zeta(s)$ satisfies a functional equation
$$
\Psi(\zeta(s), \zeta(2s), \hdots, \zeta(ds), p_1^{-s}, p_2^{-s}, \hdots, p_\ell^{-s}) 
\, = \, 0
\tag{3.10}
$$
where $p_1,\hdots,p_{\ell}$ are the prime divisors of $\vert Q \vert$
and where $\Psi \in \bold Q [X_1,\hdots,X_d, Y_1,\hdots,Y_\ell|$
is a nontrivial polynomial with rational coefficients in $d+\ell$ variables.
\par

\endproclaim

Moreover, the function $\zeta(W^{\text{{\it prof}}}(Q),s)$
extends as a multivalued analytic function
with only root singularities 
in the half--plane defined by $\operatorname{Re}(s) > 0$.
Proofs are given in Section~6, see Item~12 and Theorem~17.
\par
Furthermore, we can replace $W^{\text{{\it prof}}}(Q)$ 
by various groups having the same representation zeta function,
in particular a locally finite group $W^{\text{{\it locfin}}}(Q)$
and a finitely generated group $W^{\text{{\it fingen}}}(Q)$;
see Section~7. 
\par
By contrast, in the situation of 
a semisimple compact connected Lie group $G$,
the series $\zeta(G, s)$ has completely different properties.
Some are given in Theorem~20, which puts together
results of Weyl and Mahler, from the 1920's.

\medskip

Let us particularise the situation of Theorems 2 and 3 
to the  smallest nontrivial perfect finite group,
namely to the alternating group $A_5$ of order $60$.

\proclaim{4.~Example}
For $Q = A_5$ acting in the canonical way on
a set of $d=5$ elements, 
the representation zeta function $\zeta(s) = \zeta(W^{\text{{\it prof}}}(A_5),s)$
satisfies the functional equation
$$
\aligned
\zeta(s) \, &= \, 
60^{-1-s} \Big(
\zeta(s)^5 - 10 \zeta(s)^2\zeta(3s) - 15 \zeta(s) \zeta(2s)^2
\\
&\hskip3cm
+ 30 \zeta(2s) \zeta(3s) + 30 \zeta(s) \zeta(4s) - 36 \zeta(5s) \Big) 
\\ &+ \, 
15 \times 30^{-1-s} \times 2 \times 
\Big( \zeta(s) \zeta(2s)^2 - 2 \zeta(2s) \zeta(3s)
- \zeta(s) \zeta(4s) + 2 \zeta(5s) \Big) 
\\ &+ \, 
10 \times 20^{-1-s} \times 3 \times 
\Big( \zeta(s)^2\zeta(3s) - \zeta(2s) \zeta(3s)
-2 \zeta(s) \zeta(4s) + 2 \zeta(5s) \Big) 
\\ &+ \, 
10 \times 10^{-1-s} \left( 2 + 2^{-s} \right) 
\Big( \zeta(2s) \zeta(3s) -\zeta(5s) \Big) 
\\ &+ \, 
5 \times 5^{-1-s} \left(3 + 3^{-s} \right) 
\Big( \zeta(s) \zeta(4s) - \zeta(5s) \Big) 
\\ &+ \, 
\left( 1 + 2 \times 3^{-s} + 4^{-s} + 5^{-s} \right) \zeta(5s) .
\endaligned
\tag{3.11}
$$
Then numerical computations show that
$$
\aligned
\zeta(s) \, &\sim \,
4.186576086287 - 6.740797357 \, (s-\sigma_0)^{\frac{1}{2}}
\\
& \hskip1cm
 + 5.6535295 (s-\sigma_0) - 1.421 \, 
 (s-\sigma_0)^{\frac{3}{2}} + \cdots \endaligned
\tag{3.12}
$$
near 
$$
\sigma_0 \sim 1.1783485957546400082 .
\tag{3.13}
$$
\endproclaim

Formula (3.11) makes it easy to obtain 
the first \lq\lq few\rq\rq \ terms
of $\zeta(W^{\text{{\it prof}}}(A_5),s)$ from a computer:
$$
\aligned
\zeta(s) \, &= \, 1 
\, + \, 2 \times 3^{-s}
\, + \,          4^{-s}
\, + \,          5^{-s}
\, + \, 6 \times 15^{-s}
\, + \, 3 \times 20^{-s}
\, + \, 3 \times 25^{-s}
\\
\, &+ \, 2 \times 45^{-s}
\, + \,          60^{-s}
\, + \, 19 \times 75^{-s}
\, + \, 4   \times 90^{-s}
\, + \, 9   \times 100^{-s}
\, + \, \cdots
\endaligned
\tag{3.14}
$$
Indeed, if we agree that
$$
[ 1, 1 ], [ 3, 2 ], [ 4, 1 ], [ 5, 1 ], [ 15, 6 ], 
[ 20, 3 ], [ 25, 3 ], [ 45, 2 ], [ 60, 1 ],  [ 75, 19 ],
[ 90, 4 ], [ 100, 9 ] 
$$
is a shorthand for the $12$ first terms 
of the right--hand side of (3.14),
the nonzero terms $[n,r_n]$ of $\zeta(s)$ for $n \le 10^4$ are:
$$
\aligned
&\ssize
 [ 1, 1 ], [ 3, 2 ], [ 4, 1 ], [ 5, 1 ], [ 15, 6 ], 
[ 20, 3 ], [ 25,3 ], [ 45, 2 ], [ 60, 1 ],  [ 75, 19 ], 
\\&\ssize  
[ 90, 4 ], [ 100, 9 ], [ 125, 9 ], [ 160, 2 ], [ 180, 5], 
[ 225, 12 ], [ 240, 6 ],  [ 243, 2 ], [ 250, 2 ], [ 270, 4 ], 
\\&\ssize
  [ 300, 12 ], [ 320, 1 ], [ 375,60 ], [ 400, 3 ], [ 405, 6 ],
  [ 450, 12 ], [ 500, 28 ], [ 540, 2 ], [ 625, 27 ], [ 640, 2 ], 
\\&\ssize
[675, 2 ], [ 729, 4 ], [ 800, 6 ], [ 810, 4 ], [ 900, 52 ], 
[ 972, 2 ], [ 1024, 1 ], [ 1080, 4 ], [1125, 55 ], [ 1200, 54 ],
\\&\ssize
[ 1215, 16 ], [ 1250, 8 ], [ 1280, 4 ], [ 1350, 20 ], [ 1440, 4 ], 
[1500, 81 ], [ 1600, 12 ],  [ 1620, 12 ], [ 1875, 189 ], [ 2000, 27 ], 
\\&\ssize
[ 2025, 24 ], [ 2160, 7], [ 2250, 52 ], [ 2400, 4 ], [ 2430, 6 ],
[ 2500, 94 ], [ 2700, 30 ], [ 3000, 2 ], [ 3072, 2 ], [3125, 85 ], 
\\&\ssize
[ 3200, 6 ],  [ 3375, 18 ], [ 3600, 20 ], [ 3645, 16 ], [ 3750, 2 ], 
[ 3840, 13 ],[ 4000, 24 ], [ 4050, 40 ],  [ 4096, 1 ], [ 4320, 8 ],
\\&\ssize  
[ 4500, 339 ], [ 4800, 4 ], [ 4860, 10 ],
[ 5120, 4 ], [ 5400, 48 ],
[ 5625, 225 ], [ 5760, 4 ], [ 6000, 333 ], [ 6075, 92 ], [ 6250, 30], 
\\&\ssize
[ 6400, 18 ], [ 6480, 8 ],[ 6750, 92 ], [ 7200, 36 ], [ 7500, 442 ], 
[ 8000, 75 ], [ 8100, 106], [ 8640, 3 ], [ 9000, 12 ], [ 9375, 603 ], 
\\&\ssize
[ 9600, 20 ], [ 9720, 2 ], [ 10000, 165 ]
\endaligned
$$

\noindent
There are $2752$ non-trivial coefficients of degree $\le 10^{12}$, 
and it would be easy to extend the computations further.
\par
In Section~8, we explain our numerical computations.

\head{\bf
4.~Proof of Theorem~1, and examples involving small groups $Q$
}\endhead

The proof of Theorem~1 is 
a simple application of Arthur Clifford's theory \cite{Clif--37}
which, in particular, provides a description of the dual $\widehat G$
of a  group $G$ given as an extension
$$
1 \hskip.2cm \longrightarrow \hskip.2cm
B  \hskip.2cm \longrightarrow \hskip.2cm
G  \hskip.2cm  \longrightarrow \hskip.2cm
Q  \hskip.2cm \longrightarrow \hskip.2cm
1 
\tag{4.1}
$$
with finite quotient $Q$.
Let us recall this description in our case,
which is simpler than the general case in two respects:
the sequence splits,  $G = B \rtimes Q$,
and the representations of the subgroups of $Q$ which are involved
are linear (rather than, more generally,  projective).
Clifford's description of an irreducible representation of $G = H \wr_X Q$
has two ingredients.
\par

The first ingredient is a $Q$--orbit in $\widehat B$,
represented by some irreducible representation~$\rho$.
Since $B = H^X$, the dual $\widehat B$ can be identified with $(\widehat H)^X$,
and $\rho$ can be written uniquely 
as an outer tensor product
$\boxtimes^{x \in X} \rho_x$, 
with $\rho_x : H \longrightarrow GL(V_x)$ in $\widehat H$ for each $x \in X$.
The {\it stability subgroup} of $\rho$ is the subgroup
$$
Q_{\rho} \, = \, \left\{ q \in Q \hskip.1cm \vert \hskip.1cm
q\rho \sim \rho \right\} ,
\tag{4.2}
$$
where $\sim$ indicates equivalence of representations.
We have
$$
Q_{\rho} \, = \, Q_{P_{\rho}} ,
\tag{4.3}
$$
where the right--hand side is defined by (2.2),
and where the partition $P_{\rho}$ is that for which
$x,y \in X$ are in the same block if and only if $\rho_x \sim \rho_y$.
Moreover, since the action of $Q$ on $\widehat{B} = (\widehat H)^X$
is induced by the action of $Q$ on $X$, the representation 
$\rho = \boxtimes^{x \in X} \rho_x$ of $B$ in
the vector space $\otimes^{x \in X} V_x$
extends to a representation $\rho'$ of $B \rtimes Q_{\rho}$ 
in the same space, 
$\rho'(b,s)$ being 
for all $b \in B$ and $s \in Q_{\rho}$
the composition of $\rho(b)$ 
and the linear extension of the permutation $s$
(which commute with each other,
by definition of $Q_{\rho}$).
\par

The second ingredient is 
an irreducible representation $\sigma \in \widehat{Q_{\rho}}$.
We view it  as an irreducible representation $\sigma'$ of 
$B \rtimes Q_{\rho}$, of which the group $Q_{\rho}$ is a quotient,
so that $\rho' \otimes \sigma'$ is also 
an irreducible representation of $B \rtimes Q_{\rho}$.
Denote by
$$
\pi_{\rho,\sigma} \, = \, 
\operatorname{Ind}_{B \rtimes Q_{\rho}}^G 
\hskip.1cm (\rho' \otimes \sigma')
\tag{4.4}
$$
the induced representation, and observe that
$$
\deg \pi_{\rho,\sigma} \, = \, 
[G : B \rtimes Q_{\rho}] 
\hskip.1cm \deg \rho 
\hskip.1cm \deg \sigma   \hskip.1cm.
\tag{4.5}
$$

\proclaim{5.~Proposition (Clifford)} With the notation above,
the representation $\pi_{\rho,\sigma}$ of $G = H \wr_X Q$ is irreducible,
and any irreducible representation of $G$ is of this form, in a unique way.
\par
In other words, 
the dual $\widehat G$ is fibred (as a set) 
over the orbit space $Q \backslash \widehat B$;
the fibre over an orbit represented by $\rho \in \widehat B$ is 
the dual $\widehat{Q_{\rho}}$.
\par
In particular,  $G$ is rigid if and only if $H$ is rigid.
\endproclaim

[The general case of a group extension is more complicated, 
since $\rho'$ is a projective representation of 
$G_{\rho}$  which need not be linear,
where $G_{\rho}$ is the inverse image in $G$ of the group $Q_{\rho}$.
One has to choose $\sigma$ as a projective representation of $Q_{\rho}$
such that the class of $\rho'$ in $H^2(G_{\rho},\bold C^*)$
is minus the pull--back of the class  of  $\sigma$
in $H^2(Q_{\rho},\bold C^*)$.
In our particular case,  $G_{\rho} = B \rtimes Q_{\rho}$.
Yet there is a formulation of the proposition which
carries over to the general case of (4.1).]

\medskip

We are now ready to check the formula of Theorem~1.
\par

Consider a partition $P = (P_1,\hdots,P_j)$ in $\Pi_Q(X)$.
Let $\widehat B ^{\le}_P$ be the subset of $\widehat B$ 
consisting of those $\boxtimes^{x \in X} \rho_x$ for which $\rho_x \sim \rho_y$
as soon as $x,y$ are in the same block of $P$. 
We have a product formula
$$
\sum_{\rho = (\rho_x)_{x \in X} \in \widehat B ^{\le}_P}
\deg \left( \boxtimes^{x \in X} \rho_x \right)^{-s} 
\, = \, 
\zeta_H (\vert P_1 \vert s) \cdots \zeta_H(\vert P_j \vert s) .
$$
Let  $\widehat B ^{=}_P$ be the subset of $\widehat B^{\le}_P$
consisting of those $\boxtimes^{x \in X} \rho_x$ for which $\rho_x \sim \rho_y$
if and only if $x,y$ are in the same block of $P$.
By the defining property of the M\"obius function of $\Pi_Q(X)$,
we have
$$
\sum_{\rho = (\rho_x)_{x \in X} \in \widehat B ^{=}_P}
\deg \left( \boxtimes^{x \in X} \rho_x \right)^{-s} 
\, = \, 
\sum_{P' = (P'_1,\hdots,P'_\ell) \ge P} \mu_X(P,P')
\zeta_H (\vert P'_1 \vert s) \cdots \zeta_H(\vert P'_\ell \vert s) .
$$
It follows that the contribution to $\zeta_G (s)$
of the representations $\pi_{\rho,\sigma}$, 
with $\rho \in \widehat B^{=}_P$ and $\sigma \in \widehat{Q_P}$
(recall that $Q_P$ has been defined in (2.2)), is
$$
\sum_{\rho \in Q_P \backslash \widehat B^{=}_P, \sigma \in \widehat{Q_{\rho}}}
\deg \pi_{\rho,\sigma} ^{-s}
\, = \, 
[Q : Q_P]^{-1-s} \zeta_{Q_P} (s)
\sum_{P'  \ge P} \mu_X(P,P')
\zeta_H (\vert P'_1 \vert s) \cdots \zeta_H(\vert P'_{\ell} \vert s) ,
$$
where $P' = (P'_1,\hdots,P'_\ell)$.
The factor $[Q : Q_P]^{-s}$ is due to the induction
from $B \rtimes Q_P$ to~$G$, 
which multiplies the degrees of representations by $[Q : Q_P]$;
as we have to count only 
one $\rho$ by $Q$--orbit in $\widehat B$,
or more precisely here 
one $\rho$ by $Q_P$--orbit in $\widehat B^{=}_P$, 
there is an extra factor $[Q : Q_P]^{-1}$ 
on the right--hand side.
A summation over $P \in \Pi_Q(X)$ gives rise to the formula of Theorem~1.

\medskip

In general, computing the M\"obius function of a lattice,
for example of $\Pi_Q(X)$, is a tedious problem;
but special cases can be worked out.
For example, when the finite group $Q$ is abelian
and acts on itself ($X = Q$) by multiplications,
the computation of $\mu_X$ goes back to \cite{Dels--48}.
\par

Let us consider here some easy specialisations of Theorem~1.
For an integer $d \ge 2$, we denote by
$I_d$ the finite set $\{1,\hdots,d\}$,
by $C_d$ the cyclic group of order $d$ 
acting on $I_d$ by cyclic permutations,
by $S_d$ the symmetric group of $I_d$,
and by $H \wr_d \cdots$ the corresponding wreath products.
We denote by $\mu$ the usual M\"obius function 
of elementary number theory.

\proclaim{6.~Proposition} With the notation above we have
for an integer $d \ge 2$:
$$
\zeta(H \wr_d C_d,s) \, = \,
\sum_{e \vert d}  \left(\frac{d}{e}\right)^{-1-s} e
\sum_{f \hskip.cm \text{with} \hskip.1cm e \vert f \vert d}
\mu\left( \frac{f}{e} \right)  \zeta_H(fs) ^{d/f}                  
                                 \hskip.1cm .
\tag{4.6}
$$
In particular, if $d=p$ is prime, the summation has three terms:
$$
\zeta(H \wr_p C_p, s) \, = \,
p^{-1-s} \Big(  \zeta_H(s) ^p -  \zeta_H(ps) \Big)   
\, + \, p \zeta_H( ps)                         
                                 \hskip.1cm .
\tag{4.7}
$$
If 
$d=4$, the summation has five terms
$$
\aligned
\zeta(H \wr_4 C_4,s) \, &= \,
4^{-1-s} \Big(      
         \zeta_H(s) ^{4} - \zeta_H(2s) ^2
\Big)  
\\ \, &+ \,
2^{-s} \Big(  
        \zeta_H(2s) ^{2} -\zeta_H(4s) 
\Big)     
\, + \, 4 \zeta_H ( 4s)  
                                 \hskip.1cm .
\endaligned
\tag{4.8}
$$
For the permutation group on three objects, we have
$$
\aligned
\zeta(H \wr_3 S_3, s) \, &= \, 
6^{-1-s} \zeta_H(s)^3 +
\left( - \hskip.1cm 6^{-s}/2 + 2 \times 3^{-s} \right) 
\zeta_H(2s) \zeta_H(s)
\\
\, &+ \, 
\left(
6^{-s}/3 - 2 \times 3^{-s} + 2^{-s} + 2 \right) 
\zeta_H(3s) .
\endaligned
\tag{4.9 = 2.5}
$$
For the permutation group on four objects, 
we have
$$
\aligned
&
\zeta(H \wr_4 S_4,s) \, = \, 
\\
&\hskip.5cm 
24^{-1-s} \left(
\zeta_H(s)^4 - 6 \zeta_H(s)^2\zeta_H(2s) + 8 \zeta_H(s) \zeta_H(3s)
+ 3 \zeta_H(2s)^2 - 6 \zeta_H(4s) \right) 
\\
&\hskip.5cm 
\, + \, 
6 \times 12^{-1-s} \times 2 \times
\left( \zeta_H(s)^2 \zeta_H(2s) - 2 \zeta_H(s)\zeta_H(3s) - 
\zeta_H(2s)^2 + 2 \zeta_H(4s) \right)
\\
&\hskip.5cm 
\, + \, 
4  \times 4^{-1-s}   (2+2^{-s}) \hskip.1cm
\left( \zeta_H(s)\zeta_H(3s) - \zeta_H(4s) \right)
\\
&\hskip.5cm 
\, + \, 
3  \times 6^{-1-s} \times  4 \hskip.1cm
\left( \zeta_H(2s)^2 - \zeta_H(4s) \right)
\\
&\hskip.5cm 
\, + \, 
\left( 2 + 2^{-s} + 2 \times 3^{-s} \right) \zeta_H(4s) .
\endaligned
\tag{4.10}
$$
\endproclaim

\demo{Remark}
The reader will have no problem to guess
the formula for $\zeta(H \wr_5 A_5,s)$ from Formula (3.11);
see  (8.3).
\enddemo

\demo{Proof}
For an integer $d \ge 2$, 
the subgroups of the cyclic group $C_d$ of order $d$
are in one--to--one correspondence with the positive divisors of $d$.
If $C_d$ acts on itself by multiplication, $X = C_d$,
the orbits of a subgroup $C_e$ of $C_d$, with $e \vert d$, 
are of the form $(j, j+d/e, j+2d/e, ....)$.
If $e,f$ are positive divisors of $d$, 
the partition $P$ into orbits of the subgroup $C_e$ 
is a refinement
of the Partition $P'$ in the orbits of the subgroup $C_f$
if and only if $e \vert f$, 
or equivalently if and only if $C_e \le C_f \le C_d$.
When this is the case, 
$\mu_X(P,P') = \mu(f/e)$,
with $\mu$ the standard M\"obius function.
Thus Formula (4.6)  for $\zeta(H \wr_X C_d,s)$
is indeed a particular case of the formula of Theorem~1.
\par

If $d=p$ is prime, observe that 
the summation in (4.6) has only three terms:
$$
(d,f) \, = \,  (1,1), \hskip.2cm (1,p), \hskip.2cm (p,p) ,
$$
and (4.7) follows.
\par

If $d=4$, there are five terms:
$$
(d,f) \, = \,  (1,1), \hskip.2cm (1,2), 
\hskip.2cm (2,2), \hskip.2cm (2,4), \hskip.2cm (4,4) ,
$$
and (4.8) follows; 
the pair $(d,f) = (1,4)$ does not contribute,
because $\mu(4) = 0$.
\par

For $Q = S_3$ acting on $I_3 = \{1,2,3\}$,
the lattice $\Pi_Q(I_3)$ consists of five partitions:
the partition $\widehat 0$ in singletons, three partitions $P^{(j)}$
in the singleton $\{j\}$ and a block of size two, $j=1,2,3$,
and the partition $\widehat 1$ in one block.
The values of the M\"obius function $\mu_X$ are given
by the following table
$$
\matrix
\mu(\widehat 0,\widehat 0) &= 1 &&
   \mu(\widehat 0,P^{(j)}) &= -1 && 
       \mu(\widehat 0,\widehat 1) &= \phantom{-}2
\\
&&& 
   \mu(P^{(j)},P^{(j)}) &= \phantom{-}1 &&
       \mu(P^{(j)},\widehat 1) &=  -1
\\
&&&  
   &&&  
         \mu(\widehat 1, \widehat 1) &= \phantom{-}1 .
\endmatrix
$$
[Observe that the subgroups of $S_3$ 
of the form $(S_3)_P$
for some $P \in \Pi_{S_3}(I_3)$ are the subgroups 
of order $1$, $2$, and $6$, but not the subgroup of order $3$;
indeed, the stability subgroup of the orbit of the subgroup of order $3$
is the whole group $S_3$.]
The left--hand side of Formula (2.4) specialises to
$$
\aligned
&6^{-1-s} \Big(
\zeta_H(s)^3 - 3 \zeta_H(s) \zeta_H(2s) + 2 \zeta_H(3s) \Big)
\\
+ \hskip.2cm
&3 \times 3^{-1-s} \times 2 \Big( \zeta_H(s) \zeta_H(2s) - \zeta_H(3s) \Big)
\\
+ \hskip.2cm
&(2 + 2^{-s}) \zeta_H(3s) ,
\endaligned
$$
namely to (4.9) after minor reorganisation.
\par

Let us finally check Formula (4.10).
The lattice $\Pi_{S_4}(I_4)$ has
\roster
\item"---"
$1$ partition with  blocks of size $1,1,1,1$,
\item"---"
$6$ partitions with blocks of size $2,1,1$,
\item"---"
$4$ partitions with blocks of size $3,1$,
\item"---"
$3$ partitions with blocks of size $2,2$,
\item"---"
$1$ partition with one block of size $4$,
\endroster
namely altogether $15$ partitions.
We leave it to the reader to compute the M\"obius function.
The summation in (4.10) has $5 + 4 + 2 + 2 + 1 = 14$ terms,
more precisely:
\par

Terms with $(S_4)_P = \{1\}$,
and therefore with $\zeta_{(S_4)_P}(s) = 1$, contribute 
$$
24^{-1-s} \Big(
\zeta_H(s)^4 - 6 \zeta_H(s)^2\zeta_H(2s) + 8 \zeta_H(s) \zeta_H(3s)
+ 3 \zeta_H(2s)^2 - 6 \zeta_H(4s) \Big) .
$$
Terms with $(S_4)_P \cong S_2$ 
(fixing two of the four points of $I_4$), 
and therefore with $\zeta_{(S_4)_P}(s) = 2$, 
contribute 
$$
6 \times 12^{-1-s} \times 2 \hskip.1cm
\Big( \zeta_H(s)^2 \zeta_H(2s) - 2 \zeta_H(s)\zeta_H(3s) - 
\zeta_H(2s)^2 + 2 \zeta_H(4s) \Big).
$$
Terms with $(S_4)_P \cong S_3$,
and therefore with $\zeta_{(S_4)_P}(s) = 2 + 2^{-s}$,  
contribute 
$$
4  \times 4^{-1-s}   (2+2^{-s}) \hskip.1cm
\Big( \zeta_H(s)\zeta_H(3s) - \zeta_H(4s) \Big) .
$$
Terms with $(S_4)_P \cong S_2 \times S_2$ contribute
$$
3  \times 6^{-1-s} \times  4 \hskip.1cm
\Big( \zeta_H(2s)^2 - \zeta_H(4s) \Big) .
$$ 
The term with $(S_4)_P = S_4$ contributes
$$
\left( 2 + 2^{-s} + 2 \times 3^{-s} \right) \zeta_H(4s) .
$$
\hfill $\square$
\enddemo

\head{\bf
5.~Reminder on representations of profinite groups
}\endhead

Concerning  representations of profinite groups,
Claim (i) in the following proposition 
\footnote{
We are grateful to Bachir Bekka for showing to us
his personal notes on this.
}
is well--known, but we did not find any convenient reference.
Claim (ii) is a straightforward consequence of the definitions
(and is a particular case of, for example, 
\cite{Wils--98, Proposition 1.2.1}).

If $G$ is a profinite group 
(or more generally a compact group, for example a finite group!), 
recall that any representation of $G$ in a Hilbert space 
is unitarisable.
Recall also that, in this paper, 
representations are continuous, unless explicitly stated otherwise.

\proclaim{7.~Proposition} Let $G$ be a profinite group.
\par
(i) Let $\pi : G \longrightarrow \Cal U (\Cal H)$
be an irreducible unitary representation of $G$
in a Hilbert space~$\Cal H$.
Then $\pi$ factors through a finite quotient of~$G$.
\par
(ii)
Assume that $G = \varprojlim F_n$ is the inverse limit
of a system
$\left( p_{n,n-1} : F_n \longrightarrow F_{n-1} \right)_{n \ge 1}$
of finite groups and epimorphisms (here indexed by integers).
If $F$ is a finite group and $p : G \longrightarrow F$
a continuous epimorphism,
then there exists an integer $n$
such that $p$ is the composition of the canonical epimorphism
$p_n : G \longrightarrow F_n$ and some epimorphism
$F_n \longrightarrow F$.
\endproclaim

\demo{Proof} 
Recall that a {\it profinite group} is a compact topological group
in which every neighbourhood of $1$ contains 
an open normal subgroup
(a subgroup of a compact group which is open
is necessarily of finite index).
Equivalently, a profinite group is 
an inverse limit of finite groups.

\par
(i) We can assume that $\Cal H \ne \{0\}$.
Consider a nonzero vector $\xi \in \Cal H$.
Since $\pi$ is continuous at $1$, there exists
a neighbourhood $U$ of $1$ in $G$ such that
$$
\Vert \pi(g)\xi-\xi \Vert\, < \, 1
\hskip.5cm \text{for all} \hskip.2cm
g \in U .
$$
Since $G$ is profinite, $U$ contains
a normal subgroup $N$ of finite index.
Thus
$$
\Vert \pi(g)\xi-\xi \Vert\, < \, 1
\hskip.5cm \text{for all} \hskip.2cm
g \in N .
$$
This implies that $\eta \Doteq \int_N \pi(n) \xi dn$ 
(with $dn$ the Haar measure on $N$ of mass one)
is a nonzero $N$--invariant vector,
and in particular that the space
$\Cal H^N$ of $N$--invariant vectors
is not~$\{0\}$.
\par
Since $N$ is normal, the space $\Cal H^N$ is $\pi(G)$--invariant.
Since $\pi$ is irreducible, $\Cal H^N = \Cal H$;
in other words, $\pi$ factors through $G/N$.
\par

(ii)
Since the kernel of $p$ is closed and of finite index, it is also open. 
It follows from the definition of the topology of $\varprojlim F_n$
that $\ker p$ contains $\ker p_n$ for some $n$,
and the claim follows.
 \hfill $\square$
\enddemo

\proclaim{8.~Corollary}
Let $Q$ be a finite group acting on a finite set $X$, as in Sections~2 and 3.
We have
$$
r_n(W(Q,k+1)) \, \ge \, r_n(W(Q,k))
\tag{5.1}
$$
for all $n \ge 1$, $k \ge 0$, and
$$
r_n(W^{\text{{\it prof}}}(Q)) \, = \, 
\lim_{k \to \infty} r_n(W(Q,k)) \, \in \, \{0,1,2,\hdots,\infty\}
\tag{5.2}
$$
for all $n \ge 1$. 

Moreover, 
$r_n(W^{\text{{\it prof}}}(Q)) = 0$ unless $n$ is of the form
$p_1^{e_1} \hdots p_\ell^{e_\ell}$,
where $p_1,\hdots,p_\ell$ are the prime factors of the order of $Q$;
and the same holds for $r_n(W(Q,k))$ for all $k \ge 0$.
\endproclaim

\demo{Proof}
Inequalities (5.1) follow from (3.5) 
and Equality (5.2) follows from Proposition~7.
The last statement follows from (3.2) and from the general fact
according to which the deg\-rees of the irreducible representations
of a finite group divide the order of this group.
\newline \phantom{a}
\hfill $\square$
\enddemo

\proclaim{9.~Particular cases} 
For $d \ge 2$, any irreducible representation of the group 
$W^{\text{{\it prof}}}(S_d)$ of automorphisms of the $d$--ary tree
factors through $W(S_d,k)$, for some $k \ge 0$.
\par

Similarly, any irreducible representation 
of the group $W^{\text{{\it prof}}}(C_d)$
of $d$--adic automorphisms of the $d$--ary tree
factors through $W(C_d,k)$,
for some $k \ge 0$.
\endproclaim

In the case of $Q = C_p$ cyclic of prime order $p$,
we can be more specific.
For a finite $p$--group $H$, 
let $p^{\delta_{\operatorname{max}}(H)}$
denote the maximum of the degrees of the irreducible representations of $H$.
A first rather straightforward consequence of Formula~(4.7) is
$$
\delta_{\max}\left( W(C_p,k) \right) \, = \, 
\left\{ \aligned
1  \hskip2cm 
   &\text{if $p=2$ and $k=2$,}
\\
2^{k-2} + 2^{k-3} - 1 \hskip.5cm 
   &\text{if $p=2$ and $k \ge 3$,}
\\
1 + p + \cdots + p^{k-2} \hskip.5cm 
   &\text{if $p \ge 3$ and $k \ge 2$}.
\endaligned \right.
\tag{5.3}  
$$
A second set of consequences of (4.7), 
using slightly more calculus, is that
$$
\lim_{k \to \infty} r_{p^j}(W(C_p,k)) \, = \, \infty 
\hskip.5cm \text{for all} \hskip.2cm j \ge 0
\tag{5.4}
$$
and that the degree set of $W(C_p,k)$ is
$$
cd\left( W(C_p,k) \right) \, = \, 
\left\{p^j \hskip.1cm \vert \hskip.1cm
0 \le j \le \delta_{\max}\left( W(C_p,k) \right) \right\}
.
\tag{5.5}
$$
Formula (3.8) for $r_n(W^{\text{{\it prof}}}(C_p))$ follows.
\par

Recall that the {\it degree set} of a finite group $G$ is defined by
${cd(G) = \{n \in \bold N \hskip.1cm \vert \hskip.1cm r_n(G) > 0 \}}$.
About (5.5), let us recall that
the possible degree sets of finite $p$--groups are known;
indeed, by a theorem of Isaacs,
any finite subset of $\bold N$ of the form
$$
\{p^{e_j} \hskip.1cm \vert \hskip.1cm 0 \le j \le m\}
\hskip.5cm \text{with} \hskip.2cm 
e_0 = 0 < e_1 < \cdots < e_m
$$
is such a set \cite{Hupp--98, p.~352}.
Much less seems to be known about the degree sets
of more general finite groups; see in particular
\cite{Hupp--98, Remarks 24.5 and \S~27}.
Possible degree patterns 
(namely possible representation zeta functions of finite groups)
are even more mysterious \cite{Hupp--98, \S~6}.

\bigskip

Observe that $W^{\text{{\it prof}}}(C_d)$ 
also has {\it non--continuous} unitary representations.
For example, it can be seen that the abelianisation 
of $W^{\text{{\it prof}}}(C_d)$ is isomorphic to
the direct product $\prod_{k \ge 1} C_d$
of infinitely many copies of $C_d$. 
Given any free ultrafilter $\omega$ on $\bold N$
and a character $\chi \ne 1$ of $C_d$, 
the composition of the abelianisation 
$W^{\text{{\it prof}}}(C_d) \longrightarrow \prod_{k \ge 1} C_d$
with the $\omega$--limit
$\prod_{k \ge 1} C_d \longrightarrow \bold C^*$,
$(c_k)_{k \ge 1} \longmapsto \lim_{\omega} \chi(c_k)$
is a discontinuous character 
$W^{\text{{\it prof}}}(C_d) \longrightarrow \bold C^*$.

\head{\bf
6.~Proof of Theorems 2 and 3
}\endhead

Let $Q$ be a finite group acting on a finite set $X$ of size $d \ge 2$.
In this long section, we shall denote by $W(Q)$ the profinite group
denoted by $W^{\text{{\it prof}}}(Q)$ in Section~3. 
We  will first prove Theorem~2:
\roster
\item"---"
$W(Q)$ is rigid if and only if  $Q$ is perfect;
\endroster
and then Theorem~3, namely, in case $Q$ is perfect, that:
\roster
\item"---"
the Dirichlet series $\zeta(W(Q),s)$ 
has a finite abscissa of convergence, 
say $\sigma_0 = \sigma_0(W(Q))$,
see Proposition~15,
\item"---"
$\sigma_0 > 0$,
see Proposition~16,
\item"---"
$\zeta(W(Q),s)$ has near $\sigma_0$ a Puiseux expansion
of the form $\sum_{n=0}^{\infty} a_n (s-\sigma_0)^{n/e}$,
for some $e \le d$,
see Theorem~17.
\endroster

\proclaim{10.~Observation} 
(i) For two integers $n \ge 2$ and $d \ge 2$,
we have 
$$
\left[ \root d\of{n} \right] \, \le \, \frac{n}{2} .
\tag{6.1}
$$
\par
(ii) For a pair of integers $e,f$ such that $0 \le f \le e$
and a prime $p$, we have
$$
\frac{2+e}{(2+e-f)p^f} \hskip.2cm \left\{
\aligned
&= \, 1 \hskip.5cm \text{if} \hskip.2cm f = 0,
\\
&\le \, \frac{3}{4} \hskip.5cm \text{if} \hskip.2cm f \ge 1 .
\endaligned\right.
\tag{6.2}
$$
\endproclaim

\demo{Proof}
(i) If $n \ge 3$ and $d \ge 3$, the inequality follows from
$2^d = 8 \times 2^{d-3} < 9 \times n^{d-3} \le n^{d-1}$.
If $n \ge 4$ and $d=2$, ditto from $\sqrt{n}\le \frac{n}{2}$. 
If $n = 3$ and $d = 2$, 
then $\left[\sqrt{n} \hskip.1cm \right] = 1 < 1.5 = \frac{n}{2}$.
If $n = 2$ and $d \ge 2$, 
then $\left[ \root d\of{n} \hskip.1cm \right] 
= 1 = \frac{n}{2}$.

(ii) If $f=0$, the equality is obvious.
If $f=1$ (so that $e \ge 1$), we have
$$
\frac{2+e}{(1+e)p} \, \le \, \frac{ \frac{3}{2}(1+e)}{(1+e)p} 
\, \le \, \frac{3}{4} .
$$
If $f=2$ (so that $e \ge 2$), we have
$$
\frac{2+e}{ep^2} \, \le \, \frac{2}{p^2} \, \le \, \frac{1}{2} .
$$
If $f \ge 3$ and $f \le \frac{3}{4}e$, we have
$$
\frac{2+e}{(2+e-f)p^f} \, \le \, \frac{2+e}{(2 + \frac{1}{4}e)8} \, \le \, \frac{1}{2} .
$$
If $f \ge 3$ and $f \ge \frac{3}{4}e$, we have
$$
\frac{2+e}{(2+e-f)p^f} \, \le \, \frac{2e}{2 \times 2^{3e/4}} 
\, \le \, \frac{3}{4} .
$$
\hfill $\square$
\enddemo

For all positive $\nu \in \bold R$ and $k \in \bold N$, set
$$
r_{\nu,k} \, = \, 
\left\{\aligned
r_n(W(Q,k)) \hskip.3cm &\text{if} \hskip.2cm \nu = n \in \bold N ,
\\
0 \hskip1cm &\text{if} \hskip.2cm \nu \notin \bold N ,
\endaligned \right.
\hskip.5cm \text{and} \hskip.5cm
r_{\nu} \, = \, 
\left\{\aligned
r_n(W(Q)) \hskip.3cm &\text{if} \hskip.2cm \nu = n \in \bold N ,
\\
0 \hskip1cm &\text{if} \hskip.2cm \nu \notin \bold N .
\endaligned \right.
$$
(Values for $\nu \notin \bold N$ will only occur in Lemma~14.)
\par

In all what follows, we assume that $d = \vert X \vert \ge 2$,
and that $Q$ acts transitively on $X$.
(The transitivity hypothesis could most likely be weakened,
however arbitrary actions 
-- and in particular points of $X$ fixed by $Q$ --
would introduce unnecessary complications.)

\proclaim{11.~Lemma}
Assume that the group $Q$ is perfect and acts transitively
on a set $X$ of size $d \ge 2$.
For any $n \ge 1$, we have
$$
r_n \, = \, r_{n,\ell} \, = \, r_{n,k}
\hskip.5cm \text{for all} \hskip.2cm \ell \ge k
\hskip.2cm \text{and} \hskip.2cm  k  > \log_2 n .
\tag{6.3}
$$
\endproclaim

\demo{Proof}
{\it Step one.}
If $Q$ is perfect, then $W(Q,k)$ is also perfect,
so that $r_{1,k} = 1$ for all $k \ge 1$.
Hence, the lemma holds for $n=1$.

\medskip

{\it Step two.} 
From Clifford's theory applied to
$W(Q,k+1) = W(Q,k) \wr_X Q$,
we see that any irreducible representation $\pi$ of $W(Q,k+1)$
of degree $n$ is of the form. 
$$
\pi \, = \,
\pi_{\rho_1, \hdots,\rho_d;\sigma} \, = \, 
\operatorname{Ind}_{W(Q,k)^X \rtimes Q_P}
^{W(Q,k)^X \rtimes Q = W(Q,k+1)}
\left( \rho_1 \boxtimes \cdots \boxtimes \rho_{d} 
\boxtimes \sigma \right) ,
\tag{6.4}
$$
where $(\rho_1,\hdots,\rho_d)$ is 
a $d$--uple of irreducible representations of $W(Q,k)$, 
where $P$ is the partition of $X$ 
for which $\rho_i \sim \rho_j$ if and only if 
$i,j$ are in the same block of $P$,
and where $\sigma$ is an irreducible representation of $Q_P$
(compare with (4.4)).
Set
$$
f \, = \, 
[W(Q,k+1) : W(Q,k)^X \rtimes Q_P] \deg \sigma \, = \,  
[Q : Q_P] \deg \sigma
$$
and distinguish two cases.
\par

In the first case, $f=1$.
Thus, on the one hand $[Q : Q_P] = 1$,
and therefore $\rho_1 = \cdots = \rho_d$,
and on the other hand  $\deg \sigma = 1$, 
and therefore $\sigma = 1$. Hence
$$
n \, = \, \deg \pi_{\rho_1,\hdots,\rho_1;1} \, = \, (\deg \rho_1)^{d}  
\hskip.5cm \text{or} \hskip.5cm
\deg \rho_1 \, = \, \root d\of{\deg \pi} \, = \, \root d\of{n} .
\tag{6.5}
$$
In the second case, $f \ge 2$. Hence
$$
n \, = \,  \deg \pi_{\rho_1,\hdots,\rho_d;\sigma}  
\, \ge \, 
2 \deg(\rho_1) \cdots \deg(\rho_{d}) 
\, \ge \,
2 \deg(\rho_i) 
$$
namely
$$
\deg \rho_i  \, \le \, \frac{n}{2}
\hskip.5cm \forall i \in \{1,\hdots,d\} .
\tag{6.6}
$$
By (6.1),
$$
 \deg \rho_i \,\le \,  \left[ \frac{n}{2} \right]
 \hskip.5cm \forall i \in \{1,\hdots,d\}
$$
in the two cases.

\medskip

{\it Step three.}
Let us show  that the lemma holds for $n < 2^m$, by induction on $m$.
Since the case $m=1$ is covered by Step one,
we can assume that $m > 1$ and that the lemma holds up to $m-1$.
Step two shows that there exists a formula  of the type
$$
r_{n,\ell} \, = \, 
F \left(r_{1,\ell-1},\hdots,r_{[n/2],\ell-1} \right) 
\tag{6.7}
$$
where $F$ is an expression independent of $\ell$.
Since $[n/2] < 2^{m-1}$, we have
$$
r_{1,\ell-1} \, = \, r_{1,k-1}, \hskip.2cm 
\hdots, \hskip.2cm
r_{[n/2],\ell-1} = r_{[n/2],k-1}
$$
for $\ell-1 \ge k-1$ and $k-1 > m-1$,
by the induction hypothesis. 
It follows from (6.7) that
$$
r_{n,\ell} \, = \, r_{n,k} \hskip.5cm
\text{for all} \hskip.2cm \ell \ge k 
\hskip.2cm \text{and} \hskip.2cm k > \log_2 n
$$
whenever $n < 2^m$,
and this completes the induction step.
\hfill $\square$
\enddemo

\proclaim{12.~Proof of Theorem 2} \endproclaim

\demo\nofrills{}

If $Q$ is perfect, then $W(Q)$ is rigid by the previous lemma.
\par
Assume $Q$ is not perfect. It is known 
that the abelianisation of $W(Q,k)$
is isomorphic to the direct sum of 
$k$ copies
of the abelianisation of $Q$. 
(See \cite{BORT--96, \S~4.4, p.~145}
for a more general result, since there the action of $Q$ on $X$
need not be transitive.)
In particular, 
$$
r_1(W(Q)) \, \ge \, 
r_1\Big(\bigoplus_{j=1}^{\infty} (Q/[Q,Q])_j\Big) 
\, = \, \infty
$$
(where each $(Q/[Q,Q])_j$ denotes a copy of $Q/[Q,Q]$),
and $W(Q)$ is not rigid. 
\hfill $\square$
\enddemo

\noindent
We now proceed  to prove Theorem~3.

\proclaim{13.~Lemma}
There exists a constant $t_0 \ge 0$ with the following property.
For any pair of integers $d \ge 1$ and  $g \ge 0$, we have
$$
\sum_{g_1, \hdots, g_d \ge 0 \atop \sum_{i=1}^d g_i = g}
\left( \frac{1 + \frac{g}{2} }{ (1 + \frac{g_1}{2}) \cdots (1 + \frac{g_d}{2} ) } \right)^t
\, \le \, 4^{d-1} 
\hskip.5cm \text{for all} \hskip.2cm t \ge t_0 .
\tag{6.8}
$$
More precisely, any $t_0$ such that
$$
1 + 2 \left( \frac{5}{6} \right) ^{t_0} + \sum_{h=3}^{\infty}
\left( \frac{2}{1 + \frac{h}{2} } \right)^{t_0} \, \le \, 2
\tag{6.9}
$$
is suitable.
\endproclaim

\demo{Proof} 
For $d=1$, the lemma holds with $t_0 = 0$.
We assume from now on that $d \ge 2$, and we proceed by induction on $d$,
assuming that the lemma  holds up to $d-1$.
For $g=0$ and $g=1$, the inequality reduces respectively to
$1 \le 4^{d-1}$ and $d \le 4^{d-1}$,
so that $t_0 = 0$ is again suitable; 
we can assume therefore that $g \ge 2$.
\par

The left--hand side of (6.8) can be written as
$$
\sum_{g_1=0}^g 
\left( \frac{1+\frac{g}{2} }{(1 + \frac{g_1}{2}) (1 + \frac{g-g_1}{2})} \right)^t
\hskip.2cm
\sum_{g_2, \hdots, g_d \ge 0 \atop \sum_{i=2}^d g_i = g-g_1}
\left( \frac{1 + \frac{g-g_1}{2} }{ (1 + \frac{g_2}{2}) \cdots (1 + \frac{g_d}{2} ) } \right)^t ,
\tag{6.10}
$$
where the second summation is bounded by $4^{d-2}$, 
by the induction hypothesis.
If $g=2$, the first sum is $2 + \left( \frac{8}{9} \right)^t$, 
and in particular is bounded by $4$,
so that we can assume now that $g \ge 3$.
\par

The summation on $g_1$ from $0$ to $g$ in (6.10)
is bounded by twice the summation on $g_1$ 
from $0$ to $[g/2]$, namely by
$$
2 \left\{
1 \, + \, 
\left( \frac{1 + \frac{g}{2} }{ \frac{3}{2} ( 1 + \frac{g-1}{2} )}\right)^t
\, + \, 
\left( \frac{1 + \frac{g}{2} }{ 2 ( 1 + \frac{g-2}{2} )}\right)^t
\, + \, 
\sum_{g_1 = 3}^{[g/2]}
\left( \frac{1 + \frac{g}{2} }{ (1+\frac{g_1}{2})(1 + \frac{g-g_1}{2}) }\right)^t
\right\} .
$$
Since $g \ge 3$, we have firstly
$$
1 + \frac{g}{2} \, \le \, \frac{5}{8} + \frac{g}{8} + \frac{g}{2} 
\, = \, 
\frac{5}{6} \times \frac{3}{2} \left( \frac{1}{2} + \frac{g}{2} \right)
\hskip.5cm \text{and therefore} \hskip.5cm
\frac{1 + \frac{g}{2} }{\frac{3}{2} (1 + \frac{g-1}{2})} 
\, \le \, 
\frac{5}{6} .
$$
We have secondly
$$
\frac{1 + \frac{g}{2} }{2 ( 1 + \frac{g-2}{2} )} 
\, = \, \frac{1}{g} + \frac{1}{2} \, \le  \, \frac{5}{6} .
$$
And we have thirdly (recall that $g_1 \le g/2$)
$$
1 + \frac{g}{2} \, \le \, 1 + (g-g_1) \, < \, 
2\left(  1 + \frac{g-g_1}{2} \right)
\hskip.5cm \text{so that} \hskip.5cm
\frac{1 + \frac{g}{2} }{(1 + \frac{g_1}{2})(1 + \frac{g-g_1}{2})}
\, < \,
\frac{2}{1 + \frac{g_1}{2}} .
$$
Hence, the left--hand side of (6.8) is bounded by
$$
2  \times \left\{ 
1 +  2 \left( \frac{5}{6} \right)^t + 
\sum_{h = 3}^{\infty} \left( \frac{2}{1 + \frac{h}{2} } \right)^t
\right\}
\times 4^{d-2} ,
$$
namely by $2 \times 2 \times 4^{d-2} = 4^{d-1}$ if $t$ is large enough,
as was to be shown.
\hfill $\square$
\enddemo

Here as in Corollary~8, we denote by  
$p_1,\hdots,p_\ell$ the prime factors of $\vert Q \vert$.
Observe that they also include the prime factors 
of the degree of any irreducible representation
of one of the groups $W(Q,k)$ or $W(Q)$.

\proclaim{14.~Lemma} 
There exists a constant $t_1 \ge 0$ with the following property. 
For any $n \ge 1$, with prime decomposition
$n = p_1^{e_1}\cdots p_\ell^{e_\ell}$ ($e_1,\hdots,e_\ell \ge 0$),
we have
$$
r_{n,k} \, \le \,
\left( \frac{n}{(1+\frac{e_1}{2}) \cdots (1+ \frac{e_\ell}{2})} \right)^{t}
\hskip.5cm \text{for all} \hskip.2cm 
t \ge t_1
\hskip.2cm \text{and} \hskip.2cm
k \ge 0 .
\tag{6.11}
$$
More precisely, any $t_1$ such that 
$$
\vert Q \vert \hskip.1cm \vert X \vert ^{\vert X \vert} 
\hskip.1cm 4^{\ell(d-1)}
\vert Q \vert \hskip.1cm \left( \frac{3}{4} \right)^{t_1}  \le 1
\hskip.5cm \text{and} \hskip.2cm
t_1 \ge t_0
\tag{6.12}
$$
(with $t_0$ as in Lemma~13) is suitable.
\endproclaim

\demo{Proof} {\it Step one.}
For an integer $n \ge 1$ and $e_1,\hdots,e_\ell$ as above, set
$$
\overline{h}_n \, = \, 
\left( \frac{n}{(1+\frac{e_1}{2}) \cdots (1+ \frac{e_\ell}{2})} \right)^{t} ,
\tag{6.13}
$$
so that we wish to show that $r_{n,k} \le \overline{h}_n$
for all $k \ge 0$.
For $\nu \in \bold R$, $\nu > 0$, $\nu \notin \bold N$,
it is convenient to set $\overline{h}_{\nu} = 0$,
so that we have obviously 
$r_{\nu,k} = \overline{h}_{\nu} \hskip.2cm \forall k \ge 0$.
\par

If $k = 0$, we distinguish two cases:
if $n=1$, then $r_{1,0} = 1$ and $\overline{h}_1 = 1$ 
(the last equality for all $t$), so that $r_{n,0} = \overline{h}_n$;
if $n \ge 2$, then $r_{n,0} = 0$ and again $r_{n,0} \le \overline{h}_n$.
We assume from now on that the lemma is proven for some $k \ge 0$, 
and we will show by induction that it holds also for $k+1$.

The irreducible representations of $W(Q,k+1)$
of the first type, see (6.5) in the proof of Lemma~11, 
contribute to $r_{n,k+1}$ by
$$
r_{n,k+1}^{(i)} \, = \, 
r_{n^{1/d},k}
$$
(recall that $r_{n^{1/d},k} = 0$ if $n^{1/d}$ is not an integer).
The representations 
of the second type, see (6.6) in the proof of Lemma~11, 
contribute to $r_{n,k+1}$ by
$$
r_{n,k+1}^{(ii)} \, \le \, 
\Big( 
\sum_{P \in \Pi_Q(X)}
\hskip.2cm \sum_{\sigma \in \widehat{Q_P}}
\Big) '
\hskip.2cm \sum_{
n_1,\hdots,n_d \ge 1, \atop
 \prod_{i=1}^d n_i = n / ((\deg \sigma)[Q:Q_P])
}
r_{n_1,k} \cdots r_{n_d,k}
$$
where the prime in $\left(\sum_{P} \sum_{\sigma}\right)'$
indicates that the pair $(P,\sigma)=(\hat 1,1)$ does not occur.
Observe that $\deg \sigma$ is a divisor of $\vert Q_P \vert$,
so that $ [Q : Q_P] \deg \sigma$ is a divisor of $\vert Q \vert$.
If we introduce the constant
$$
K \, = \, 
\text{number of pairs $(P,\sigma)$,
with $P \in \Pi_Q(X)$,  
$\sigma \in \widehat{Q_P}$,
and $(P,\sigma) \ne (\hat 1,1)$} ,
$$
it follows that
$$
r_{n,k+1} 
\, = \, 
r_{n,k+1}^{(i)} + r_{n,k+1}^{(ii)}
\, \le \, 
r_{n^{1/d},k} 
\, + \,
K  \sum_{2 \le f \le \vert Q \vert \atop f \vert n} \hskip.1cm
\sum_{n_1,\hdots,n_d \ge 1 \atop \prod_{i=1}^d n_i = n / f}
r_{n_1,k} \cdots r_{n_d,k} .
\tag{6.14}
$$
Observe that the number of choices for $P$
is strictly bounded by $\vert X \vert ^{\vert X \vert}$
and that the number of choices for $\sigma$ 
is bounded by  $\vert Q \vert$,
so that $K+1 \le \vert Q \vert \vert X \vert^{\vert X \vert}$.
\medskip

{\it Step two.}
In the last sum of Inequality $(6.14)$, any $n_i$ which occurs,
namely any $n_i$ such that $r_{n_i,k} \ne 0$,
is a product of the $p_j$'s. 
We repeat the definition of the exponents $e_j$,
and we define exponents $e_{i,j}$, $f_j$ by
$$
n \,  = \, \prod_{1 \le j \le \ell} p_j^{e_j} , 
\hskip.5cm
n_i \, = \, \prod_{1 \le j \le \ell} p_j^{e_{i,j}} 
\hskip.2cm \text{for} \hskip.2cm i=1,\hdots,d, 
\hskip.5cm
f    \, = \, \prod_{1 \le j \le \ell} p_j^{f_j} .
\tag{6.15}
$$
We use the induction hypothesis to bound 
the first term of (6.14):
$$
r_{n^{1/d},k} \, \le \, 
\overline{h}_{n^{1/d}} 
\hskip.1cm \cdots \hskip.1cm
\overline{h}_{n^{1/d}} 
\hskip.1cm
\overline{h}_{n^{1/d}/p_j} ,
$$
where we have $d-1$ factors $\overline{h}_{n^{1/d}}$
and one factor in which $j$ is such that $p_j$ divides $n$
(in case there does not exist any such $j$,
we have $r_{n^{1/d},k} = 0$).
We use again the induction hypothesis to bound 
the second term of (6.14),
and we collect terms to obtain
$$
r_{n,k+1} \, \le \,
(K+1)
\sum_{2 \le f \le \vert Q \vert \atop f \vert n} \hskip.1cm
\sum_{n_1,\hdots,n_d \ge 1 \atop f \prod_{i=1}^d n_i = n }
\overline{h}_{n_1} \cdots \overline{h}_{n_d} .
$$
Using the definition of the $\overline{h}_{n_i}$
and reordering the terms, we have
$$
\aligned
r_{n,k+1} \, &\le \,
(K+1)
\sum_{2 \le f \le \vert Q \vert \atop f \vert n} \hskip.1cm
\sum_{n_1,\hdots,n_d \ge 1 \atop f \prod_{i=1}^d n_i = n }
\\
& \hskip3cm
\left( \frac{n_1}{ (1 + \frac{e_{1,1}}{2}) \cdots (1 + \frac{e_{1,\ell}}{2} ) } \right)^t
\cdots
\left( \frac{n_d}{ (1 + \frac{e_{d,1}}{2}) \cdots (1 + \frac{e_{d,\ell}}{2} ) } \right)^t
\\
&= \,
(K+1) 
\sum_{2 \le f \le \vert Q \vert \atop f \vert n} \hskip.1cm
\sum_{n_1,\hdots,n_d \ge 1 \atop f \prod_{i=1}^d n_i = n } \hskip.1cm
\prod_{i=1}^d \Big( n_i \Big)^t  \hskip.1cm
\prod_{j=1}^\ell \left( \frac{1}{ 1 + \frac{e_j}{2} } \right)^t
\\
& \hskip3cm
\left( \frac{ 1 + \frac{e_1}{2} }
{(1 + \frac{e_{1,1}}{2}) \cdots (1 + \frac{e_{d,1}}{2})} \right)^t
\cdots
\left( \frac{ 1 + \frac{e_\ell}{2} }
{(1 + \frac{e_{1,\ell}}{2}) \cdots (1 + \frac{e_{d,\ell}}{2})} \right)^t 
\\
&= \,
(K+1) \left( \frac{n}{ \prod_{j=1}^\ell (1 + \frac{e_j}{2}) }  \right)^t
\sum_{2 \le f \le \vert Q \vert \atop f \vert n} 
\frac{1}{f^t}  \hskip.1cm
\sum_{n_1,\hdots,n_d \ge 1 \atop f \prod_{i=1}^d n_i = n } 
\\
& \hskip3cm
\left( \frac{ 1 + \frac{e_1}{2} }
{(1 + \frac{e_{1,1}}{2}) \cdots (1 + \frac{e_{d,1}}{2})} \right)^t
\cdots
\left( \frac{ 1 + \frac{e_\ell}{2} }
{(1 + \frac{e_{1,\ell}}{2}) \cdots (1 + \frac{e_{d,\ell}}{2})} \right)^t
\endaligned
$$
where the term
$ \left( \frac{n}{ \prod_{j=1}^\ell (1 + \frac{e_j}{2}) }  \right)^t$
is precisely $\overline{h}_n$.
We replace now a sum of products by a product of sums,
and we obtain
$$
\aligned
r_{n,k+1} \, &\le \,
 (K+1) \hskip.1cm \overline{h}_n
\sum_{2 \le f \le \vert Q \vert \atop f \vert n} 
\frac{1}{f^t}  \hskip.1cm \times \hskip.1cm
\left(
\sum_{e_{1,1},\hdots,e_{d,1} \ge 0 \atop \sum_{i=1}^d e_{i,1} = e_1 - f_1}
\left( 
\frac{ 1 + \frac{e_1}{2} }{ (1 + \frac{e_{1,1}}{2} ) \cdots (1 + \frac{e_{d,1}}{2}) } 
\right) ^t
\right)
\times
\\
& \hskip4cm
\cdots \hskip.2cm \times
\left(
\sum_{e_{1,\ell},\hdots,e_{d,\ell} \ge 0 \atop \sum_{i=1}^d e_{i,\ell} = e_\ell - f_\ell}
\left( 
\frac{ 1 + \frac{e_\ell}{2} }{ (1 + \frac{e_{1,\ell}}{2} ) \cdots (1 + \frac{e_{d,\ell}}{2}) } 
\right) ^t
\right)
\\
&= \,
 (K+1) \hskip.1cm \overline{h}_n
\sum_{2 \le f \le \vert Q \vert \atop f \vert n} 
\left(
\frac{ 1 + \frac{e_1}{2} }{ \left( 1 + \frac{e_1 - f_1}{2} \right) p_1^{f_1} }
\right)^t
\hskip.1cm \cdots \hskip.1cm 
\left(
\frac{ 1 + \frac{e_\ell}{2} }{ \left( 1 + \frac{e_\ell - f_\ell}{2} \right) p_\ell^{f_\ell} }
\right)^t
\hskip.1cm \times
\\
& \hskip3cm
\left(
\sum_{e_{1,1},\hdots,e_{d,1} \ge 0 \atop \sum_{i=1}^d e_{i,1} = e_1 - f_1}
\left( 
\frac{ 1 + \frac{e_1-f_1}{2} }{ (1 + \frac{e_{1,1}}{2} ) \cdots (1 + \frac{e_{d,1}}{2}) } 
\right) ^t
\right)
\times
\\
& \hskip4cm
\cdots \hskip.2cm \times
\left(
\sum_{e_{1,\ell},\hdots,e_{d,\ell} \ge 0 \atop \sum_{i=1}^d e_{i,\ell} = e_\ell - f_\ell}
\left( 
\frac{ 1 + \frac{e_\ell-f_\ell}{2} }{ (1 + \frac{e_{1,\ell}}{2} ) \cdots (1 + \frac{e_{d,\ell}}{2}) } 
\right) ^t
\right) .
\endaligned
$$
Each of the $\ell$ terms
$$
\frac{ 1 + \frac{e_1}{2} }{ \left( 1 + \frac{e_1 - f_1}{2} \right) p_1^{f_1} }
\hskip.2cm , \hskip.5cm
\hdots , \hskip.2cm
\frac{ 1 + \frac{e_\ell}{2} }{ \left( 1 + \frac{e_\ell - f_\ell}{2} \right) p_\ell^{f_\ell} }
$$
above is bounded by $1$,
and at least one of them is bounded by $\frac{3}{4}$ 
(see Observation 10.ii).
Each of the next $\ell$ sums over $d$--uples of $e_{\star,\star}$'s 
is bounded by $4^{d-1}$, by Lemma~13.
It follows that
$$
\aligned
r_{n,k+1} \, &\le \, 
\overline{h}_n \hskip.1cm (K + 1) \hskip.1cm 4^{\ell(d-1)}
\sum_{2 \le f \le \vert Q \vert \atop f \vert n} \left( \frac{3}{4} \right)^t
\\
\, &\le \, 
\overline{h}_n \hskip.1cm (K + 1) \hskip.1cm 4^{\ell(d-1)}
\vert Q \vert \hskip.1cm \left( \frac{3}{4} \right)^t .
\endaligned
$$
For $t$ large enough, this shows that 
$$
r_{n,k+1} \, \le \, \overline{h}_{n} 
$$
and ends the induction argument.
\hfill $\square$
\enddemo

\proclaim{15.~Proposition}
Let $Q$ be a perfect finite group acting transitively 
on a finite set $X$ with at least two points. 
Then the representation zeta function $\zeta(W(Q), s)$
has a finite abscissa of convergence, say $\sigma_0(W(Q))$.
\par
Moreover, the function $\zeta(W(Q), s)$, 
which is holomorphic in the half--plane
defined by $\operatorname{Re}(s) > \sigma_0(W(Q))$,
has a singularity at $\sigma_0(W(Q))$.
\endproclaim

\demo{Proof}
It is elementary to check that
$r(W(Q)) = \sum_{n=1}^{\infty} r_n(W(Q)) = \infty$
(if necessary, see the first step of the proof of the next proposition,
which is independent of the present proof).
It follows that $\sigma_0(W(Q))$ is given by Formula (1.2).
Hence, by Lemma 14, and using the notation 
$\overline{h}_n$ of (6.13):
$$
\aligned
\sigma_0(W(Q)) 
\, &\le \, \limsup_{n \to \infty} \frac{ \ln \left( \sum_{j=1}^n \overline{h}_j \right)}{\ln n}
\\
\, &\le \, \limsup_{n \to \infty} \frac{1}{\ln n}
\ln \left( \sum_{d_1=0}^{e_1} \cdots \sum_{d_\ell = 0}^{e_\ell}
\left( \frac{n}{ (1 + \frac{d_1}{2}) \cdots (1 + \frac{d_\ell}{2}) }  \right)^t \right)
\\
\, &\le \, \limsup_{n \to \infty} \frac{1}{\ln n}
\ln \left( \sum_{d \hskip.1cm \vert \hskip.1cm n} n^t \right)
\, \le \, 
\limsup_{n \to \infty} \frac{1}{\ln n} \ln \left( n \times n^t \right)
\\
\, &\le \, t+1 
\endaligned
$$ 
where $e_1,\hdots,e_\ell$ are defined in terms of $n$ as in Lemma~14.
Observe that, at this stage, the only clear bounds are
$0 \le \sigma_0(W(Q)) < \infty$.
\hfill $\square$
\enddemo

\proclaim{16.~Proposition} 
In the situation of the previous proposition
($Q$ perfect acting transitively on $X$, 
with $d = \vert X \vert \ge 2$), we have
$$
\sigma_0(W(Q)) \, > \, 0 .
\tag{6.16}
$$
\endproclaim

\demo{Proof}
{\it Step one: there exist infinitely many values of $n$ such that
$r_n(W(Q)) \ge 1$.}
\par

Indeed, since $Q = W(Q,1)$ is perfect and not reduced to $\{1\}$,
there exists $n \ge 2$ such that $r_{n}(Q) \ge 1$.
Let $k \ge 1$  be such that $r_n(W(Q,k)) \ge 1$
and let $\rho_1$ be 
an irreducible representation of $W(Q,k)$ of dimension $n$.
The $n^d$--dimensional representation 
$\pi_{\rho_1,\hdots,\rho_1;1}$, as in (6.4),
contributes to $r_{n^d}(W(Q,k+1)) \ge 1$.
Continue with $n^d$ in lieu of $n$.
\smallskip

{\it Step two: there exist infinitely many values of $n$ 
for which $r_n(W(Q)) \ge 2$.}
\par

Let $N \ge 2$ be an integer.
By Step one, there exist $k \ge 1$ and irreducible representations
$\rho_1, \hdots, \rho_d$ of $W(Q,k)$ of pairwise distinct degrees, all at least $N$.
Since $Q$ is perfect, 
and therefore strictly contained in the symmetric group of $X$,
the action of $Q$ is not $d$ times transitive.
Hence there exists a permutation $\tau$ of $X$ such that
$$
\rho_1, \hdots, \rho_d
\hskip.5cm \text{and} \hskip.5cm
\rho_{\tau(1)}, \hdots, \rho_{\tau(d)}
$$
are {\it not} in the same orbit of $Q$ 
acting on $\left(\widehat{W(Q,k)}\right)^d$.
It follows that the irreducible representations 
$$
\operatorname{Ind}_{W(Q,k)^X}^{W(Q,k+1)}
\left( \rho_1 \boxtimes \cdots \boxtimes \rho_d \right)
\hskip.5cm \text{and} \hskip.5cm
\operatorname{Ind}_{W(Q,k)^X}^{W(Q,k+1)}
\left( \rho_{\tau(1)} \boxtimes \cdots \boxtimes \rho_{\tau(d)} \right)
$$
of $W(Q,k+1)$, which are both of degree 
$\vert Q \vert \prod_{i=1}^d \deg \rho_i$ (a degree $> N^d$),
are not equivalent.
The claim of Step two follows.
\smallskip

{\it Step three: For any positive integer $B$, 
there exist infinitely many values of $n$
for which $r_n(W(Q)) \ge B$.}
\par

We proceed by induction on $B$ (see Step two for $B=2$).
Suppose that Step three has been shown for some value $B_0 \ge 2$.
Let $N \ge 2$ be an integer.
By the induction hypothesis, 
there exist integers $k \ge 1$  and $n_1,\hdots,n_d$
such that $N < n_1 < \cdots < n_d$ and $r_{n_i}(W(Q,k)) \ge B_0$, 
$i=1,\hdots,d$.
For each $i$, choose $r_{n_i}(W(Q,k))$ 
pairwise inequivalent irreducible representations
$\rho_{i,j_i}$ of $W(Q,k)$ of degree $n_i$.
The irreducible representations
$$
\operatorname{Ind}_{W(Q,k)^X}^{W(Q,k+1)} \left(
\rho_{1,j_1} \boxtimes \cdots \boxtimes \rho_{d,j_d} 
\right) ,
\hskip1cm 
1 \le j_i \le r_{n_i}(W(Q,k)) \hskip.5cm (i=1,\hdots d)
$$
of $W(Q,k+1)$ are pairwise inequivalent,
all of the same degree, which is $\vert Q \vert \prod_{i=1}^d n_i > N^d$,
and there are $\prod_{i=1}^d r_{n_i}(W(Q,k)) \ge (B_0)^d$ of them.
\smallskip

{\it Step four: Set $K = \root d\of{ \frac{1}{d! \, 2^d} }$.
Let $n \ge 1$ be such that
$r_n(W(Q)) \ge \min \{ 2d, \frac{2}{K} \}$.
Then
$$
r_{\vert Q \vert n^d} (W(Q)) \, \ge \, \Big( K \,  r_n(W(Q)) \Big) ^d .
$$
}
\par

Observe that the existence of the integer $n$ 
involved in the claim of Step four follows from Step three.
Let $k \ge 1$ be such that there exist 
pairwise inequivalent irreducible representations
$\rho_1, \hdots, \rho_{r(n)}$ of $W(Q,k)$ of dimension $n$,
where we have written $r(n)$ for $r_n(W(Q))$.
Any choice of $d$ distinct representations 
$\rho_{j_1}, \hdots, \rho_{j_d}$
among $\rho_1, \hdots, \rho_{r(n)}$ 
provides an irreducible representation
$$
\operatorname{Ind}_{W(Q,k)^X}^{W(Q,k+1)}
\left( \rho_{j_1} \boxtimes \cdots \boxtimes \rho_{j_d} \right)
$$
of $W(Q,k+1)$ of dimension $\vert Q \vert n^d$.
Different choices provide nonequivalent representations,
and there are
$\binom{r(n)}{d} \ge \frac{(r(n)/2)^d}{d!} = ( K r(n) )^d$
such choices.
\smallskip

{\it Step five: end of proof.}
\par

Choose $B$ such that $\frac{d}{d-1} \ln K + \ln B > 0$.
By iteration of the inequality shown in Step four, for some $m \ge 2$, we have
$$
\aligned
& r(m) \ge B , \hskip.5cm
r(\vert Q \vert m^d) \ge (K r(m))^d , \hskip.5cm 
\\
& \hskip1cm
r(\vert Q \vert^{1+d} m^{d^2}) \ge \left( K (K r(m))^d \right)^d 
= K^{d+d^2} r(m)^{d^2}, \hskip.5cm \hdots
\endaligned
$$
and more generally
$$
r\left( \vert Q \vert^{(d^t - 1)/(d-1)} m^{d^t} \right) 
\, \ge \,
K^{d(d^t-1)/(d-1)} \,  B^{d^t}
\hskip.5cm \text{for any integer} \hskip.2cm t \ge 1 .
$$
Hence
$$
\aligned
\sigma_0(W(Q)) \, &= \, 
\limsup_{n \to \infty} \frac{ \ln \left( \sum_{j=1}^n r(j) \right) }{ \ln n}
\\
&\ge \,
\limsup_{t \to \infty}
\frac{ 
\ln \Big( r\big( \vert Q \vert^{(d^t - 1)/(d-1)} m^{d^t} \big) \Big)
}{ 
\ln \Big( \vert Q \vert^{(d^t - 1)/(d-1)} m^{d^t} \Big)
} 
\\
&\ge \,
\limsup_{t \to \infty}
\frac{ 
\ln \Big( K^{d(d^t-1)/(d-1)} B^{d^t} \Big)
}{ 
\ln \Big( \vert Q \vert^{(d^t - 1)/(d-1)} m^{d^t} \Big)
} 
\\
&= \,
\limsup_{t \to \infty}
\frac{ 
d\frac{d^t-1}{d-1} \ln K + d^t \ln B
}{
\frac{d^t-1}{d-1} \ln \vert Q \vert + d^t \ln m
}
\\
&= \,
\frac{
\frac{d}{d-1}  \ln K + \ln B
}{
\frac{1}{d-1} \ln \vert Q \vert + \ln m
}
\endaligned
$$
and the last fraction is positive by the choice of $B$.
\hfill $\square$
\enddemo

We concentrate briefly on the group $Q=A_5$ to
introduce the next proposition.
Consider the polynomial  
$\Psi \in \bold Q [X_1,X_2,X_3,X_4,X_5,Y_2,Y_3,Y_5]$
defined by
$$
\aligned
&
\Psi(X_1,X_2,X_3,X_4,X_5,Y_2,Y_3,Y_5) \, = \,
\\
&\hskip.5cm 
\phantom{\, = \,} 
\frac{1}{60} Y_2^2 Y_3 Y_5
\Big(
X_1^5 - 10 X_1^2 X_3 - 15 X_1 X_2^2
+ 30 X_2 X_3 + 30 X_1 X_4 - 36 X_5 
\Big) 
\\
&\hskip.5cm 
\, + \, 
Y_2Y_3Y_5
\Big( 
X_1 X_2^2 -2 X_2X_3 - X_1X_4 + 2 X_5
\Big) 
\\
&\hskip.5cm 
\, + \, 
\frac{3}{2} Y_2^2Y_5
\Big( 
X_1^2 X_3 - X_2X_3 -2 X_1X_4 +2 X_5
\Big) 
\\
&\hskip.5cm 
\, + \, 
Y_2Y_5 \left( 2 + Y_2 \right) 
\Big( X_2 X_3 - X_5 \Big) 
\\
&\hskip.5cm 
\, + \, 
Y_5 \left(3 + Y_3 \right) 
\Big( X_1 X_4 - X_5 \Big)
\\
&\hskip.5cm 
\, + \, 
\left( 1 + 2   Y_3 + Y_2^2 + Y_5 \right) X_5 
\\
&\hskip.5cm 
\, - \, X_1.
\endaligned
\tag{6.17}
$$
Since $W(A_5) \wr_5 A_5 \cong W(A_5)$,
the representation zeta function $\zeta(s) = \zeta(W(A_5),s)$
is a solution of the equation
$$
\Psi \left(
\zeta(s), \zeta(2s), \zeta(3s), \zeta(4s), \zeta(5s), 2^{-s}, 3^{-s}, 5^{-s} 
\right) 
\, = \, 0
\tag{6.18}
$$
(compare with (3.11)).
\par

Let $X_\ast$ stand for $\{X_2,X_3,X_4,X_5\}$ and
$Y_\ast$ for $\{Y_2,Y_3,Y_5\}$; we can also write
$$
\Psi(X_1,X_\ast,Y_\ast) \, = \, \sum_{i=0}^5 v_i(X_\ast,Y_\ast)X_1^{5-i} ,
$$
with
$$
v_0(X_\ast,Y_\ast) \, = \, \frac{1}{60} Y_2^2Y_3Y_5 ,
\hskip.5cm 
v_3(X_\ast,Y_\ast) \, = \, 
- \frac{1}{6} X_3 Y_2^2Y_3Y_5  + \frac{3}{2} X_3 Y_2^2 Y_5 ,
\hskip.5cm \hdots
$$
(it happens that $v_1(X_\ast,Y_\ast) = v_2(X_\ast,Y_\ast) = 0$,
but this does not play any role in the argument below).
Let $R$ denote the ring of holomorphic functions 
in the half--plane of inequation 
$\operatorname{Re}(s) > \frac{1}{2} \sigma_0$,
where $\sigma_0   \Doteq \sigma_0(W(A_5))$.
Since the Dirichlet series for $\zeta(2s), \hdots, \zeta(5s)$ 
converge when $\operatorname{Re}(s) > \frac{1}{2} \sigma_0$,
we can view $v_0,\hdots,v_5$ as elements of $R$ and
$\Psi(...)$ as a one variable polynomial $\Phi(X_1)$ in the ring $R [X_1]$,
so that $\zeta(s)$ is now a root of the polynomial $\Phi(X_1)$.
Near $\sigma_0$,  all coefficients of this polynomial are bounded;
moreover,  the top coefficient is
an entire function which is never zero;
if follows that $\zeta(s)$ tends 
to a finite limit $\zeta(\sigma_0)$ when
$s \to \sigma_0$ (say by real values $s$ with $s > \sigma_0$).
Thus $\zeta(s)$ is a root of a polynomial of the form
$$
\Theta(X) \, = \, X^5 + \sum_{i=1}^5 u_i(s) X^{d-i}
\tag{6.19}
$$
where the coefficients $u_i$ are in the ring $R$.
Let 
$$
\Delta \, \Doteq \, 
\operatorname{Disc}(\Theta(X)) \, = \, 
\prod_{1 \le i < j \le 5} (\alpha_i - \alpha_j)^2
\, \in \, R 
\tag{6.20}
$$
denote the  {\it discriminant} of $\Theta$,
where $\alpha_1,\hdots,\alpha_5$ are the roots of $\Theta$
in an appropriate extension
of the field of fractions of $R$
(for discriminants, see e.g.\ \cite{Bour--81, 
chapitre IV, \S~6, no~7, proposition~11}).
Since $\Delta$ is a function  
which is holomorphic in a neighbourhood of $\sigma_0$, 
we can evaluate $\Delta$ at $\sigma_0$
to find its value $\Delta(\sigma_0) \in \bold C$.

\medskip

We return to the general case. 
We write  $\zeta(s)$ for the representation zeta function $\zeta(W(Q),s)$,
and $\sigma_0$ for its abscissa of convergence $\sigma_0(W(Q))$.
For $k \ge 0$, denote by $\Cal H(k)$ the half--plane defined by
$\operatorname{Re}(s) > 2^{-k}\sigma_0$ and by
$\Cal O (\Cal H(k))$ the ring of holomorphic functions there.
Observe that
the Dirichlet series defining $\zeta(s)$ converges in $\Cal H (0)$,
and that 
$$
\Cal H (0) \, \subset \, 
\Cal H (1)  \, \subset \,  \cdots  \, \subset \, 
\Cal H (k)  \, \subset \,  \hdots   \, \subset \, 
\bigcup_{k=0}^{\infty} \Cal H(k) \, = \, 
\left\{ 
s \in \bold C \hskip.1cm \vert \hskip.1cm \operatorname{Re}(s) > 0
\right\} .
$$
\par

\proclaim{17.~Theorem}
Let $Q$ be a perfect group acting transitively on a set $X$ with
$d\ge 2$ points, as in Propositions 15 and 16,
and let $\zeta(s)$, $\sigma_0$, $\Cal H (1)$ be  as above.
Recall that $p_1,\hdots,p_{\ell}$ denote the prime divisors of the order of $Q$.
\par

(i) There exists a polynomial 
$\Psi \ne 0$ in $\bold Q [X_1,X_2,\hdots,X_d,Y_1,Y_2,\hdots,Y_\ell]$ such that
$$
\Psi(\zeta(s), \zeta(2s), \hdots, \zeta(ds), p_1^{-s}, p_2^{-s}, \hdots, p_\ell^{-s}) 
\, = \, 0 .
\tag{6.21 = 3.10}
$$
\par

(ii) There exists a polynomial 
$$
\Theta(X) \, = \, X^d + \sum_{i=1}^d u_i(s) X^{d-i} \in \Cal O (\Cal H (1)) [X]
\tag{6.22}
$$
with discriminant $ \Delta(s) \in \Cal O (\Cal H (1))$
such that $\zeta(s)$ is a root of $\Theta$ and such that
$\Delta(\sigma_0) = 0$.
\par

(iii) When $s \to \sigma_0$ with $\operatorname{Re}(s) > \sigma_0$,
the limit $\lim_{s \to \sigma_0} \zeta(s)$ exists and is finite.
\par

(iv) Near $\sigma_0$, 
the function $\zeta(s)$ has a Puiseux expansion of the form
$$
\zeta(s) \, = \, \sum_{n=0}^{\infty} a_n (s-\sigma_0)^{n/e} ,
\tag{6.23}
$$
where $e$ is an integer such that $2 \le e \le d$.
\par

(v)
The function $\zeta(s)$ extends as a multivalued analytic function
with only root singularities 
in the half--plane defined by $\operatorname{Re}(s) > 0$.
\endproclaim

\demo{Proof}
(i)
In view of the isomorphism (3.7),
this follows from Theorem~1.

(ii)
For the definition of $\Theta$ and for $\zeta(s)$ being a root of $\Theta$,
see the argument for the case of $Q = A_5$ just before the theorem.
If $\sigma_0$ were not a zero of $\Delta$, 
we could extend holomorphically $\zeta(s)$ to the left of~$\sigma_0$, 
in contradiction with the definition of~$\sigma_0$ 
and the classical result, recalled in the introduction just before (1.2),
that $\sigma_0$ is a singular point of $\zeta(s)$.
\par

(iii)
As $\zeta(s)$ near $\sigma_0$ satisfies a monic polynomial equation
with all coefficients bounded,
it has a finite limit when $s \to \sigma_0$ 
with $\operatorname{Re}(s) > \sigma_0$.
\par

(iv)
This is standard in the Newton--Puiseux theory~;
see for example \cite{Casa--00, Theorems 1.5.4 and 1.7.2}.
Observe that one could not have $e=1$ in (6.23), 
since, then, the point $\sigma_0$ would not be 
a singularity of~$\zeta(s)$.
\par

(v) We can define inductively for each $k \ge 0$
a multivalued analytic extension of $\zeta(s)$ in $\Cal H (k)$,
such that $\zeta(s)$ is a root of a polynomial
formally identical to (6.22),
but now with coefficients analytic in $\Cal H (k+1)$.
Claim (v) follows.
\hfill $\square$
\enddemo

The argument of (ii) shows that, moreover,
$\zeta(s)$ has a continuous extension to the closed half--plane
$\operatorname{Re}(s) \ge \sigma_0$.

\medskip
\noindent
The proof of Theorem~3 is now complete.

\medskip

\proclaim{18.~Question}
Is  it true that,
for any curve in the half--plane $\{\operatorname{Re}(s) > 0 \}$
avoiding the singularities and ending on the imaginary axis,
the function $\zeta(s)$ of Theorem~17
does {\rm not} extend holomorphically
along this curve up to the end point?

In other words, do all continuations of $\zeta(s)$ admit
$\{\operatorname{Re}(s) = 0 \}$ as a natural frontier?
\endproclaim

\proclaim{Numerical estimates}
When  $Q = A_5$ acts on $X = I_5$,
a numerical estimate for $\sigma_0 = \sigma_0(W(Q))$ 
is given in (3.13).
The first terms 
\footnote{
We are grateful to Don Zagier for having computed these terms for us.
}
of the Puiseux expansion of $\zeta(W(A_5),s)$
near $\sigma_0$ are as in (3.12); in particular, $e=2$.
\endproclaim

Thus, in the case of the group $A_5$, as before Theorem~17,
we have $e = 2 < d = 5$.
It follows that the polynomial $\Theta(X)$,
now regarded with coefficients  in the ring of germs
of holomorphic functions around $\sigma_0$,
is reducible, with an irreducible factor of degree $e=2$.

A few numerical experiments suggest a positive answer to Question 18. 
In the case of $W(A_5)$, it seems that the discriminant 
$\Delta$ of Claim (ii) in Theorem~17
vanishes at $\sigma_0$ and two other places
$$
s'_0 \pm  it'_0 \,  \sim \, 0.8973038819 \pm 0.0264098303 \, i 
\tag{6.24}
$$
in $\{s \in \bold C \hskip.1cm \vert \hskip.1cm 
\operatorname{Re}(s) > \sigma_0/2\}$.
\par

A small program was written to do the following: a small complex step
$\epsilon$ is chosen, for instance $0.001+0.0001i$, and an integer $N$
large enough so that $\operatorname{Re}(N\epsilon)>2$ is chosen.
The values of $\zeta(W(A_5),s)$ at 
$s=\epsilon, 2\epsilon, \hdots, N\epsilon$ 
will be approximated by numerical values $z_n$. 
These are computed in decreasing order $z_N,z_{N-1},\hdots,z_1$. 
If $\operatorname{Re}(n\epsilon)>2$, then $z_n$ is computed 
using the first $10^{12}$ terms of the power series. Otherwise, 
the program ensures that $z_{2n}, z_{3n}, z_{4n}, z_{5n}$ are computed, 
and obtains a polynomial functional equation for $z_n$, with five roots; 
the root closest to $z_{n+1}$ is chosen as approximation for $z_n$.

If one plots in 3-space the points 
$(\operatorname{Re}(\epsilon)n,\operatorname{Re}(z_n),
\operatorname{Im}(z_n))$, 
one sees discrete approximations of continuous curves, 
that approximate $\zeta(W(A_5),s)$ 
along radial half--lines $\bold R_+\epsilon$. 
The following remarks can be made empirically:
\roster
\item"---"
these curves remain bounded as 
$\operatorname{Re}(\epsilon) \to 0$;
\item"---"
they are smooth except when $\epsilon$ is a real multiple 
of $\sigma_0$ or $\sigma'_0 \pm it'_0$;
\item"---"
for $\epsilon$ a real multiple of $1+i$, 
one sees as $\operatorname{Re}(\epsilon)n \to 0$ 
a spiral with non-vanishing radius and faster and faster winding;
\item"---"
other random values of $\epsilon$ indicate either convergence, or
oscillation, or spiralling.
\endroster

\medskip

\noindent {\bf 19.~Remark.}
Let $G$ be a connected compact Lie group 
which is semisimple and simply connected; 
denote by $\ell$ its rank 
and by $\kappa = \frac{1}{2}(\operatorname{dim}(G) - \ell)$ 
the number of positive roots
of the Lie algebra $\operatorname{Lie}(G) \otimes_{\bold R} \bold C$.
By Hermann Weyl's theory, there exists a polynomial $P$
of degree $\kappa$ in $\ell$ variables such that
$$
\zeta(G,s) \, = \, \sum_{n_1,\hdots,n_{\ell} \ge 0}
P(n_1,\hdots,n_{\ell})^{-s} .
$$
For example, if $G$ is simple of rank $\ell = 2$, we have
\cite{Hump--72, Page~140}
$$
\aligned
\zeta(SU(3),s) \, &= \,
\sum_{m_1 = 0}^{\infty} \sum_{m_2 = 0}^{\infty}
\Big( \frac{1}{2} (m_1+1)(m_2+1)(m_1+m_2+2)  \Big)^{-s} 
\\
\, &= \, 1 + 2 \cdot 3^{-s} + 2 \cdot 6^{-s} + 8^{-s} + 2 \cdot 10^{-s}
 + 4 \cdot 15^{-s} + 2 \cdot 21^{-s} + 2 \cdot 24^{-s}  
 \\
 \, &+ 27^{-s} + 2 \cdot 28^{-s} + 2 \cdot 35^{-s} + 2 \cdot 36^{-s}
  + 2 \cdot 42^{-s} + 2 \cdot 45^{-s} + 2 \cdot 48^{-s}
   \\
 \, &+ 2 \cdot 55^{-s} + 2 \cdot 60^{-s} + 2 \cdot 63^{-s} +
 64^{-s}   + 2 \cdot 66^{-s} + \cdots\cdots\cdots
\\
\zeta(Spin_5(\bold C),s) \, &= \,
\sum_{m_1 = 0}^{\infty} \sum_{m_2 = 0}^{\infty}
\Big( \frac{1}{6} (m_1+1)(m_2+1)(m_1+m_2+2)
(2m_1+m_2+3) \Big)^{-s} 
\\
\zeta(G_2,s) \, &= \,
\sum_{m_1 = 0}^{\infty} \sum_{m_2 = 0}^{\infty}
\Big(  \frac{1}{120}
(m_1+1)(m_2+1)(m_1+m_2+2) \hskip.2cm \times
\\
& \hskip2cm
(m_1+2m_2+3)(m_1+3m_2+4)(2m_1+3m_2+5)  \Big)^{-s}.
\endaligned
$$

\proclaim{20.~Theorem (Weyl and Mahler)}
Let $G$ be a connected compact Lie group 
which is semisimple and simply connected,
and let $\ell,\kappa$ be as above.

The abscissa of convergence of $\zeta(G,s)$ is $\frac{\ell}{\kappa}$
and $\zeta(G,s)$ extends to a meromorphic function in the whole
complex plane. Moreover the poles of $\zeta(G,s)$ are all simple
and at rational points of the real axis; in particular,
the largest singularity
$\sigma_0(G) = \frac{\ell}{\kappa}$ is a simple pole of $\zeta(G,s)$.
\endproclaim

\demo{References}
This theorem 
\footnote{
We are grateful to Olivier Mathieu who first made one of us
aware of this.
}
is a straightforward consequence of
several classical results. More precisely,
given the formula of Weyl for the dimensions of the irreducible
representations of $G$,
it becomes a particular case of Satz~II in \cite{Mahl--28}, 
where the more general case of
$\sum_{m_1,\hdots,m_{\ell} \ge 0}
P_1(m_1,\hdots,m_{\ell}) P_2(m_1,\hdots,m_{\ell})^{-s}$
is considered.
\par
See also  \cite{LaLu--08, Theorem 5.1},
as well as other results on  Dirichlet series of related form
in \cite{Sarg--84}.
\hfill $\square$
\enddemo

A comparison with Theorem~17 shows how different
representation zeta functions can be from each other.

\head{\bf
7.~Locally finite groups and finitely generated groups
acting on rooted trees
}\endhead

Besides the group $W^{\text{{\it prof}}}(Q)$ defined in (3.6),
we have also for each $k \ge 1$ a natural injection
$$
W(Q,k) \, \longrightarrow \, W(Q,k+1) 
\tag{7.1}
$$
(which is the natural splitting of the epimorhism (3.5)),
and therefore a corresponding locally finite group
$$
W^{\text{{\it locfin}}}(Q) \, = \, \bigcup_{k \ge 1} W(Q,k) .
\tag{7.2}
$$
\par

For integers $d \ge 2$ and $k \ge 0$, 
let $T_d(k)$ denote the {\it $d$--ary rooted tree of height $k$;}
its vertices are the finite words in $I_d = \{1,2,\hdots,d\}$ 
of length at most $k$,
and any word $w$ of length between $0$ (the empty word) and $k-1$
is adjacent to the $k$ words $wx$, for $x \in \{1,\hdots,d\}$.
We denote by $T_d$ the {\it infinite $d$--ary rooted tree},
in which the subtree induced by 
the vertices at distance at most $k$ from the root
is precisely $T_d(k)$.
\par

The particular case of (3.1) and (3.6) in which $Q$ is
the symmetric group $S_d$ acting in the standard way on 
the set $I_d \Doteq \{1,2,\hdots,d\}$  provides 
{\it the full automorphism groups of these trees:}
$$
W(S_d,k) \, = \, \operatorname{Aut}\left( T_d(k) \right)
\hskip.5cm \text{and} \hskip.5cm
W^{\text{{\it prof}}}(S_d) \, = \, \operatorname{Aut}\left( T_d \right) .
\tag{7.3}
$$
The case of $Q = C_d$ acting by cyclic permutations on $I_d$
gives rise to the so--called 
{\it group of $d$--adic automorphisms} of $T_d$.
Observe that, for $p$ a prime, 
$W(C_p,k)$ is a $p$--Sylow subgroup of $W(S_p,k)$,
and $W(C_p)$ is a $p$--Sylow subgroup of $W(S_p)$.
\medskip

{\it
For simplicity, we assume from now on that $Q$ 
acts faithfully on the set $X$ of size $d \ge 2$,
}
identified to $I_d$.
The groups defined by (3.1) and (3.6) are therefore naturally
groups of tree automorphisms:
$$
W(Q,k)  \,  \subset \,  \operatorname{Aut}\left( T_d(k) \right)
\hskip.5cm \text{and} \hskip.5cm
W^{\text{{\it prof}}}(Q)  \,  \subset \,  \operatorname{Aut}\left( T_d \right) .
\tag{7.4}
$$
The epimorphism $W(Q,k+1) \longrightarrow W(Q,k)$ of (3.5)
is the restriction to $T_d(k)$ of automorphisms of  $T_d(k+1)$,
and the epimorphism $W(Q,k) \longrightarrow Q$
corresponding in an analogous way to (3.1)
is the restriction to $X$ of  automorphisms of $T_d(k)$,
where $X$ is identified with the set of leaves of $T_d(1)$.
\par

Moreover, corresponding to the splitting which appears in (3.5),
we may extend \lq\lq rigidly\rq\rq \ automorphisms of $T_d(k)$
to automorphisms of $T_d(\ell)$, $\ell > k$, and of~$T_d$;
otherwise written:
$$
W(Q,k) \, \subset \, W(Q,k+1) \,  \subset \,  \cdots  \,  
\subset \,  W^{\text{{\it locfin}}}(Q) 
\,  \subset \, \operatorname{Aut}\left( T_d \right) .
\tag{7.5}
$$
It follows that the group defined by (7.2) is also a group
$$
W^{\text{{\it locfin}}}(Q) \, \subset \, \operatorname{Aut}\left( T_d \right)
\tag{7.6}
$$
of tree automorphisms.
\par

There are other interesting dense
subgroups of $W^{\text{{\it prof}}}(Q)$,
and in particular finitely generated groups
which play an important role in various questions 
(see e.g.\ \cite{Bart--03b}, \cite{Neum--86}, and \cite{Wils--04})
as the examples below show.
\par

It follows from Theorem~5 in~\cite{Sega--01} that, if $Q$ is perfect,
acts transitively on $X$, and distinct points have distinct
stabilizers, then $W^{\text{{\it prof}}}(Q)$ contains dense finitely
generated subgroups.
\par

{\it Assume furthermore that $Q$ is $2$--transitive on $X$,}
and consider two elements of $X$, written $1$ and $2$ in (7.9)  below.
Let 
$$
W^{\text{{\it fingen}}}(Q) \, = \, \langle Q , \overline{Q} \rangle 
\, \subset \, \operatorname{Aut}\left( T_d \right)
\tag{7.7} 
$$
be the group of automorphisms of $T_d$
generated by two copies of $Q$.
An element $q$ in the first copy $Q$ 
acts on $T_d$ as follows,
with $(x_1,\hdots,x_k)$ a typical vertex of $T_d$:
$$
q(x_1,x_2,x_3,\hdots,x_k) \, = \, (q(x_1),x_2,x_3,\hdots,x_k)
\tag{7.8}
$$
for all $k \ge 0$ and $x_1,\hdots,x_k \in X$.
An element $\overline{q}$ of the second copy $\overline{Q}$ 
acts on $T_d$ by
$$ 
\overline{q}(x_1,x_2,\hdots,x_k) \, = \, \left\{
\aligned
&
(x_1,\hdots,x_j,q(x_{j+1}),x_{j+2}, \hdots, x_k)
\\
&\hskip2.3cm
\text{if} \hskip.3cm x_1 = \cdots = x_{j-1} = 1,
\hskip.1cm x_j = 2,
\\
&
(x_1,\hdots,x_k) \hskip.5cm \text{otherwise,}
\endaligned \right.
\tag{7.9}
$$
for all $k$ and $x_1,\hdots,x_k$ as above.
\par

{\it 
We now assume furthermore that  $Q$ is generated by 
$\bigcup_{x \neq  y \in X}
\left(\operatorname{Stab}_{Q}(x) \cap 
\operatorname{Stab}_{Q}(y)\right)'$,
}
where the prime indicates a commutator subgroup;
this holds for example if $Q$ is simple and 
$\operatorname{Stab}_{Q}(x) \cap \operatorname{Stab}_{Q}(y)$ is not Abelian.
The case $Q = A_6$ and $X = I_6$ 
was considered by Peter Neumann in \cite{Neum--86}.
We define a finitely generated group
$$
W^{PN}(Q) \, = \, \langle 
\omega(x,q) \hskip.1cm \vert \hskip.1cm 
x \in X, \hskip.1cm q \in \operatorname{Stab}_{Q}(x) 
\rangle
\, \subset \, \operatorname{Aut}(T_d)
\tag{7.10}
$$
where the generators $\omega(x,q)$ are defined by
$$ 
\omega(x,q) (x_1,x_2,\hdots,x_k) \, = \, \left\{
\aligned
&
(x_1,\hdots,x_j,q(x_{j+1}),x_{j+2}, \hdots, x_k)
\\
&\hskip2.3cm
\text{if} \hskip.3cm x_1 = \cdots = x_{j} = x,
\hskip.1cm x_{j+1} \ne x,
\\
&
(x_1,\hdots,x_k) \hskip.5cm \text{otherwise.}
\endaligned \right.
\tag{7.11}
$$
\par

\proclaim{21.~Proposition}
Let $Q$ be a perfect finite group acting faithfully and transitively
on a finite set $X$ of size  $d \ge 2$.
Let $W$ be one of the groups
$$
W^{\text{{\it prof}}}(Q) , \hskip.5cm
W^{\text{{\it locfin}}}(Q) , \hskip.5cm
W^{\text{{\it fingen}}}(Q) , \hskip.5cm
W^{PN}(Q) 
$$
(with hypotheses on the action of $Q$ on $X$ as above
for the groups $W^{\text{{\it fingen}}}(Q)$ 
and $W^{\text{{\it PN}}}(Q)$).
\par

Then $W$ is perfect, residually finite, and
isomorphic to its own permutational wreath product with $Q$:
$$
W \cong W \wr_X Q.
$$
\endproclaim

\demo{Proof}
These groups are perfect because they are generated by perfect
subgroups, and are residually finite 
because they are all  subgroups  
of the profinite group $W^{\text{{\it prof}}}(Q)$.
\par

   It remains to establish the isomorphism; it comes from the natural
inclusion of $W^{\text{{\it prof}}}(Q)$ as the first factor in the base of
$W^{\text{{\it prof}}}(Q) = W^{\text{{\it prof}}}(Q) \wr_X Q$.
The isomorphism is clear for the first two examples. 
For $W^{\text{{\it fingen}}}(Q)$, see~\cite{Bart--03a, Proposition~3.7}
or~\cite{Bart--03b, Proposition~2.2}. 
For $W^{PN}(Q)$, we argue as follows.
\par

Consider 
$q,r \in  \operatorname{Stab}_{Q}(x) \cap \operatorname{Stab}_{Q}(y)$.
Then
$[\omega(x,q),\omega(y,r)]$ acts as $[q,r]$ on $T_d$ as in~(7.8). 
By our hypothesis, $W(Q)^{PN}$ contains the whole of $Q$ acting as in~(7.8). 
It then follows that $W(Q)^{PN}$ contains
$\omega(x,q)q^{-1}$ (where $q^{-1}$ acts as in~(7.8)), 
which maps to $(1,\dots,\omega(x,q),\dots,1)1$
through the embedding $W^{PN}(Q) \to W^{PN}(Q) \wr_X Q$. 
We conclude that this embedding is an isomorphism. 
This is in essence the argument
given by Neumann in \cite{Neum--86,   pp.\ 307 sqq}.
\hfill $\square$
\enddemo

\proclaim{22.~Proposition (P.~Neumann)} 
Let $G$ be a perfect, residually finite group isomorphic to $G \wr_X  Q$.
Then $G$ is just infinite.
\par
More precisely, any homomorphism $G\to H$ with non-trivial kernel
factors as $G\to W(Q,k)\to H$ for some $k\in\bold N$.
\endproclaim

\demo{Reference}
This is part of \cite{Neum--86, Theorem 5.1}.
The homomorphism $G \longrightarrow W(Q,k)$ is obtained as
$$
\aligned
G \, &\longrightarrow \, G \wr_X Q 
\, \longrightarrow \, (G \wr_X Q) \wr_X Q) = G \wr_{X^2} W(Q,2)
\, \longrightarrow \, \cdots 
\\
\,& \longrightarrow \, G \wr_{X^k} W(Q,k)
\, \longrightarrow \, W(Q,k) .
\endaligned
$$
\hfill $\square$
\enddemo

\proclaim{23.~Proposition} In the situation of Proposition~22,
there exists for any integer $n \ge 1$ an integer $k(n) \ge 0$
with the following property:
\par
Any finite dimensional representation $G \longrightarrow GL_n(\bold C)$
factors through $W(Q,k(n))$; 
in  particular $r_n(G) < \infty$ for all $n\ge 1$. 
More precisely,
$$
r_n(G) \, = \, 
\max_{k \in \bold N} r_n(W(Q,k)) \, < \, \infty.
$$
\endproclaim

\demo{Proof}
Let $\rho : G \to GL_n(\bold C)$ be a representation. 
In view of Proposition~22, 
it suffices to show that $\rho$ is not faithful.
\par

Since $G \cong G \wr_X Q$, the group $G$ contains 
for all $\ell \in \bold N$
a  finite subgroup  isomorphic to~$Q^{X^\ell}$, 
and {\it a fortiori} a subgroup isomorphic to 
$(\bold Z / p\bold Z)^{\ell}$,
where $p$ is a prime which divides the order of $G$.
By the representation theory of finite abelian groups,
the image of $(\bold Z / p\bold Z)^{\ell}$ by $\rho$
is up to conjugation a diagonal $p$--subgroup of $GL_n(\bold C)$,
so that the order of this image is at most $p^n$;
it follows that  $\rho$ has a non--trivial kernel $N$.
\par

The proof of Proposition~22,
see  \cite{Neum--86, Lemma~5.2}, actually shows that,
if $N$ is  a normal subgroup of $G$ 
with non-trivial image in $W(Q,k)$, 
then $N$ contains the kernel of the map $G \to W(Q,k)$.
\hfill $\square$
\enddemo

\demo{Remark} 
This gives another proof of 
(5.2) in Corollary~8
for $W^{\text{{\it prof}}}(Q)$.
\enddemo

From the previous propositions, we deduce finally:

\proclaim{24.~Theorem}
In the situation of Proposition~21, the groups
$$
W^{\text{{\it prof}}}(Q) , \hskip.5cm
W^{\text{{\it locfin}}}(Q) , \hskip.5cm
W^{\text{{\it fingen}}}(Q) , \hskip.5cm
W^{PN}(Q) 
$$
are rigid and have the same representation zeta function $\zeta(W^{\text{{\it prof}}}(Q),s)$,
that which is the object of Theorem~2, Theorem~3, and Section~6.
\endproclaim

\head{\bf
8.~Comments on numerical computations
}\endhead

The way chosen to write Formula (2.4) in Theorem~1
shows clearly the inclusion--exclusion ingredient of its proof.
However, it is not optimal for numerical computations,
because evaluations of $\zeta_H(\cdots)$
are much more expensive to compute that any other term,
so that the number of occurrences of these $\zeta_H(\cdots)$
should be kept minimal, as in the second way of writing below.
$$
\aligned
\zeta(H \wr_X Q,s) \, &= \,
\sum_{P \in \Pi_Q(X)} [Q : Q_P ]^{-1-s}
\zeta_{Q_P}(s) 
\\
& \hskip1.5cm
\sum_{P' = (P'_1,\hdots,P'_\ell) \ge P} \hskip.1cm \mu_X( P,P')  \hskip.1cm 
\zeta_H(\vert P'_1 \vert s) \cdots \zeta_H(\vert P'_\ell \vert s) 
\\
&= \,
\sum_{P' = (P'_1,\hdots,P'_\ell) \in \Pi_Q(X)} 
\zeta_H(\vert P'_1 \vert s) \cdots \zeta_H(\vert P'_\ell \vert s)
\\
& \hskip1.5cm
\sum_{P \le P'} \hskip.1cm \mu_X( P,P')  \hskip.1cm 
[Q : Q_P ]^{-1-s} \zeta_{Q_P}(s) .
\endaligned
\tag{8.1 = 2.4}
$$

We used GAP~\cite{GAP4} to produce specialisations of Formula~(8.1),
namely to produce Formulas~(4.6) to (4.10), 
as well as (8.3), (8.4), and (8.6) below. 
For small groups, they also have or could have been produced by hand, 
by enumerating all subgroups $S$ of $Q$; 
computing the partition $P$ of $X$ they induce; 
computing the stabilizer $Q_P$ of $P$; 
and keeping those pairs $(S,P)$ for which $S=Q_P$.

It turned out to be much faster, for the examples we considered, 
to enumerate all partitions $P$ of $X$; 
to compute their stabilizer $Q_P$; 
to compute the partition $P'$ of $X$ induced by $Q_P$; 
and to keep those $(Q_P,P)$ for which $P =P'$.

All these commands (enumerating partitions, computing stabilizers,
comparing groups) are simple instructions in GAP.

It is then straightforward, using the offspring of~(8.1), 
to compute the zeta function $\zeta(W(Q,k),s)$ for small $k$, 
as a polynomial in $p_1^{-s},\dots,p_\ell^{-s}$. 
\medskip

We describe now how numerical estimates 
as those in (3.12) and (3.13) were obtained.
We concentrate specifically on $Q=A_5$ in this section; 
the same method works {\it mutatis mutandis} for other examples.

We stress that our goal is not to prove formal enclosures for the
numerical constants $\sigma_0, \zeta(W(Q),\sigma_0), \dots$, 
but rather to obtain good approximations of their decimal expansion, 
so as to check
(e.g.\ in Plouffe's Inverter \cite{P-INV})
whether these constants already appeared in mathematics. 
We were unsuccessful with the examples we considered.

First, a large number $N$, say $10^{15}$, is chosen; 
and the degree--$(\le N)$ truncations 
$P_k(s)$ of $\zeta(W(Q,k),s)$ are computed. 
By Corollary~8, the sequence $\left(P_k(s)\right)_{k \ge 0}$ 
is eventually constant, with limit $P(s)$; 
and $P(s)$ is the degree--$(\le N)$ truncation of $\zeta(W(Q),s)$.

We noted experimentally, 
by varying $N$ between $10^{10}$ and $10^{15}$, 
that $\zeta(W(Q),s)$ is very well approximated 
by $P(s)$ when $s>\sigma_0+1$; 
for example, again for $Q = A_5$, 
the value at $s = 2.1$ of $P(s)$ 
differs only at its $a$th digit 
when comparing $N=10^a$ and $N=10^{a+1}$.

This is supported by the following heuristic. 
We know that the series $\zeta(W(Q),s)$ 
converges at $\sigma_0+\epsilon$ for all $\epsilon>0$, 
so its general term $r_n n^{-\sigma_0-\epsilon}$ is bounded; 
say, for simplicity, bounded by $1$. 
We also assume, somewhat crudely, 
that the tail $\sum_{n>N}r_n n^{-\sigma_0 - \epsilon}$ 
is bounded by $1$. 
Therefore, for $s > \sigma_0+1$, 
the error $\zeta(W(Q),s) - P(s)$  is bounded by
$$
\sum_{n>N} r_n n^{-\sigma_0-\epsilon-1} \, < \, 
N^{-1} \sum_{n>N} n^{-\sigma_0-\epsilon} \precsim N^{-1}.
$$ 
For our choice of $N$, 
we may therefore expect $P(s)$ to approximate $\zeta(W(Q),s)$
accurately to $15$ digits.

We will soon see that $\sigma_0 > 1$, 
and we will only evaluate the Dirichlet polynomial $P(s)$ at $es$ 
for $s \sim \sigma_0$ and $e \ge 2$. 
All our estimations will then be accurate to about $15$ digits.

Our first goal is to estimate $\sigma_0(W(Q))$. 
For this, in (6.17) and as in~(6.18), we replace $X_2, \dots, X_5$ 
by $P(2s), \dots, P(5s)$  respectively, 
and we replace each $Y_j$ by $j^{-s}$, 
yielding a degree--$d$ polynomial $\Phi$ in $X = X_1$ 
whose coefficients are Dirichlet polynomials in $s$.

It is computationally impractical 
to evaluate algebraically the discriminant $\Delta(s)$, 
so we resort to another trick: 
we first find an enclosure $[l_0,u_0]$ for $\sigma_0$, 
such that the number of real roots of $\Phi(X)$
differs at $s=l_0$ and at $s=u_0$. 
We then repeatedly compute the number of real roots of $\Phi(X)$
at $m_i \Doteq (u_i+l_i)/2$, 
and set $[l_{i+1},u_{i+1}] \Doteq [l_i,m_i]$ 
if the number of real roots at $m_i$ equals that at $u_i$, 
and set $[l_{i+1},u_{i+1}] \Doteq [m_i,u_i]$ otherwise. 
The enclosure's width is halved at each step, 
so after $50$ or so iterations, for our choice of $N$, 
we have found a good approximation of $\sigma_0$.

Note that, at $s=\sigma_0$, the polynomial $\Phi(X)$
has a multiple root; let $t+1$ be its multiplicity. 
We assume for simplicity that $t=1$,
which is the only case we have encountered in our numerical experiments.

We next set $s = \sigma_0$ in $\Phi(X)$, 
and obtain $a_0$, the first term 
in the Puiseux series of $\zeta(W(Q),s)$ near $\sigma_0$, see (6.23). 
We then compute for $i = 0,1,\hdots$ the 
coefficient $a_{i+1}$ as the average of 
$(\Phi(X) - a_0 - \cdots - 
a_i(s-\sigma_0)^{i/2})/(s-\sigma_0)^{(i+1)/2}$ 
over $s$ very near $\sigma_0$.

We have obtained in this manner a Puiseux expansion
$$
\zeta(W(Q),s) \, \sim \, a_0 + a_1 \sqrt{s-\sigma_0} 
+ a_2(s-\sigma_0) + \cdots
\tag{8.2}
$$
for $\zeta(W(Q),s)$ near $\sigma_0$; 
see (3.12) and (3.13) for the actual numbers,
for the case of $A_5$ acting on $I_5$.
Let us now consider a few other examples.

\bigskip 

\centerline{* * * * * * *}

\bigskip

The group $A_5$ of order $60$ is the smallest nontrivial perfect finite group.
As above, we denote a  wreath product with respect to 
its natural action on $\{1,2,3,4,5\}$  by $H \wr_5 A_5$.
We can also view $A_5 \cong PSL_2(\bold F_5)$ 
as acting on the projective line over the Galois field of order five,
and we use the notation $\wr_6$ for this wreath product.
\par

The next smallest finite simple group is $PGL_3(\bold F_2)$,
which is equal to $GL_3(\bold F_2)$ and has order $168$.
It has a natural action on the projective plane with $7$ points.
We denote a  wreath product with respect to this action 
by $H \wr_7 PGL_3(\bold F_2) $.

Let us write down the particular cases of Formula (8.1)
corresponding to these actions, and some results
of numerical computations.

\proclaim{25.~Proposition}
With the notation described above, we have
$$
\aligned
&\zeta(H \wr_5 A_5,s) \, = \, 
60^{-1-s}  \zeta_H(s)^5 
\\
&\hskip.5cm + \, 
\Big( 20^{-1-s} \times 30 - 60^{-1-s} \times 10 \Big) \zeta_H(s)^2 \zeta_H(3s)
\\
&\hskip.5cm + \, 
\Big( 30^{-1-s} \times 30 - 60^{-1-s} \times 15 \Big) \zeta_H(s) \zeta_H(2s)^2
\\
&\hskip.5cm + \, 
\Big( 5^{-1-s} \left(3 + 3^{-s} \right) 5 - 20^{-1-s} \times 60 
- 30^{-1-s} \times 30 + 60^{-1-s} \times 30 \Big) 
\zeta_H(s) \zeta_H(4s)
\\
&\hskip.5cm + \, 
\Big( 10^{-1-s} \left(2 + 2^{-s} \right) 10 - 20^{-1-s} \times 30 
- 30^{-1-s} \times 60 + 60^{-1-s} \times 30 \Big) 
\zeta_H(2s) \zeta_H(3s)
\\
&\hskip.5cm + \, 
\Big( 1 + 2 \times 3^{-s} + 4^{-s} + 5^{-s}
- 5^{-1-s}\left(3 + 3^{-s} \right) 5
- 10^{-1-s} \left( 2 + 2^{-s} \right) 10
\\
& \hskip5cm
+ 20^{-1-s} \times 60
+ 30^{-1-s} \times 60
- 60^{-1-s} \times 36 \Big) \zeta_H(5s)
\endaligned
\tag{8.3}
$$
and
$$
\aligned
\zeta(H \wr_6 PSL_2(\bold F_5),s) \, &= \, 60^{-1-s} \zeta_H(s)^6
\\
&+ \, \Big(30^{-s}-\frac{60^{-s}}{4}\Big)
                    \zeta_H(s)^2\zeta_H(2s)^2
\\
&+ \, \Big(2\cdot6^{-s}+2\cdot12^{-s}-2\cdot30^{-s}+\frac25\cdot60^{-s}\Big)
                    \zeta_H(s)\zeta_H(5s)
\\
&+ \, \Big(\frac43\cdot15^{-s}-30^{-s}+\frac{60^{-s}}{6}\Big)
                       \zeta_H(2s)^3
\\
&+ \, \Big(2\cdot10^{-s}+20^{-s}-2\cdot30^{-s}+\frac{60^{-s}}{3}\Big)
                       \zeta_H(3s)^2
\\
&+ \, \Big(1 + 2\cdot3^{-s} + 4^{-s} + 5^{-s}-2\cdot6^{-s} - 2\cdot10^{-s}
\\
& \hskip1cm - 2\cdot12^{-s}
         -\frac43\cdot15^{-s}-20^{-s}+4\cdot30^{-s}-\frac23\cdot60^{-s}\Big)
                     \zeta_H(6s) .
\endaligned
\tag{8.4}
$$
For comparison with (3.13), the Dirichlet series 
$\zeta(W(PSL_2(\bold F_5)),s)$ converges for
$$
\operatorname{Re}(s) \, > \,
\sigma_0(W(PSL_2(\bold F_5))) \, \sim \,  1.13333324(7) ,
\tag{8.5}
$$ 
again with a Puiseux expansion in $\sqrt{s-\sigma_0}$.

\medskip

We have next
$$
\aligned 
& \zeta(H \wr_7 PGL_3(\bold F_2),s) \, = \,
\left(\frac{168^{-s}}{168} \right)     \zeta_H(s)^7
\\
&\hskip.5cm + \, 
\left(\frac{84^{-s}}{2}-\frac{168^{-s}}{8}\right)
     \zeta_H(s)^3 \zeta_H(2s)^2
\\
&\hskip.5cm + \, 
\left(\frac{2}{3} 42^{-s}-\frac{84^{-s}}{2}+\frac{168^{-s}}{12}\right)
      \zeta_H(s)^3 \zeta_H(4s)
\\ 
&\hskip.5cm + \, 
\left(\frac{2}{3} 42^{-s}-\frac{84^{-s}}{2}+\frac{168^{-s}}{12}\right)
       \zeta_H(s) \zeta_H(2s)^3
\\
&\hskip.5cm + \, 
\left(4 \times 21^{-s}-3 \times 42^{-s}+\frac{84^{-s}}{2}\right)
      \zeta_H(s) \zeta_H(2s) \zeta_H(4s)
\\
&\hskip.5cm + \, 
\left(2 \times 28^{-s} + 56^{-s} - 2 \times 84^{-s} + \frac{168^{-s}}{3}\right)
      \cdot\zeta_H(s) \zeta_H(3s)^2
\\
&\hskip.5cm + \, 
\left(2 \times 7^{-s}+14^{-s} - 2 \times 21^{-s} - 2 \times 28^{-s} +
      \frac{42^{-s}}{3} - 56^{-s} + 2 \times 84^{-s} - \frac{168^{-s}}{3} \right)
\\
&\hskip.5cm
\qquad\qquad (\zeta_H(s) \zeta_H(6s)+\zeta_H(3s) \zeta_H(4s))
\\
&\hskip.5cm + \, 
\Big(1+2 \times 3^{-s} + 6^{-s} -  3 \times 7^{-s} + 8^{-s}
      -2 \times 14^{-s} + 2 \times 28^{-s} + 42^{-s} + 56^{-s} 
\\
&\hskip.5cm
\qquad\qquad   -2 \times 84^{-s} + \frac{2}{7} \times 168^{-s}\Big)
      \zeta_H(7s).
\endaligned
\tag{8.6}
$$
Moreover  $\zeta(W(PGL_3(\bold F_2)),s)$ converges for
$$
\operatorname{Re}(s) \, > \,
\sigma_0(W(PGL_3(\bold F_2))) \, \sim \,   1.112156628 ,
\tag{8.7}
$$ 
again with a Puiseux expansion in $\sqrt{s-\sigma_0}$.
\endproclaim

\head{\bf
9.~Some special values of representation zeta functions
}\endhead

Let us record a few general facts about representation zeta functions.
If $G$ is a rigid group, observe that 
$$
\zeta(G,0) \, = \, \sum_{n \ge 1} r_n(G) \, = \,  
\vert \widehat G \vert
\, \Doteq \, r(G)
\tag{9.1}
$$
is the number (possibly infinite) of equivalence classes 
of irreducible representations of $G$.
\par

These functions are well adapted to direct products:
if $G$ and $H$ are two rigid groups, 
$G \times H$ is also rigid and we have
$$
\zeta(G \times H , s) \, = \, \zeta(G,s) \zeta(H,s) .
\tag{9.2}
$$
\par

When $G$ is finite, $r(G)$ is the {\it class number} of $G$
(= its number of conjugacy classes).
In this case, observe also that the number 
of representations of dimension $1$ is
$$
\lim_{s \to \infty, s > 0} \zeta(G,s) \, = \, 
r_1(G) \, = \,  \vert G / [G,G] \vert ,
\tag{9.3}
$$
and that
$$
\zeta(G,-2) \, = \, \vert G \vert .
\tag{9.4}
$$
If $G$ is a finite group,
$\zeta(G,s)$ is an entire function:
if $G$ is a finite $p$--group, 
$\zeta(G,s)$ is a polynomial in $p^{-s}$.
\par

Let $G$ be a finite group such that
$\frac{1}{\vert G \vert} 
\sum_{g \in G} \operatorname{trace}(\pi(g^2)) = 1$
for all $\pi \in \widehat{G}$, namely such that
all $\pi \in \widehat{G}$ can be realised over $\bold R$; then we have
$$
\zeta(G,-1) \, = \,
\sharp \left\{ s \in G \hskip.1cm \vert \hskip.1cm s^2 = 1 \right\} .
\tag{9.5}
$$
See for example \cite{Serr--98, Section~13, Exercise~3}
or \cite{BeZh--98, Chapter~4, Theorems~13 and~14}.
Groups for which (9.5) holds include 
the dihedral groups $D_n$ for which 
$$
\zeta(D_{2n+1},-1) \, = \, 2n+2
\hskip.5cm \text{and} \hskip.5cm
\zeta(D_{2n},-1) \, = \, 2n+2 ,
$$ 
the symmetric groups $S_n$ (see e.g.\ \cite{CuRe--62, \S~28}),
the alternating groups $A_n$ when $n \in \{5,6,10,14\}$
(and for no other $n \ge 3$), 
the projective linear groups $PSL_2(\bold F_q)$
when $q$ is either a power of $2$ or of the form $4k+1$;
see 
\footnote{
We are grateful to Alexander Zaleskii for information about Schur indices.
}
\cite{Feit--83, Theorem 6.1}.
Observe that (9.5) cannot hold if $G$ is of non--trivial odd order.

\medskip

If we set $r(k) \Doteq r\left( W(C_2,k) \right)$, 
the  sequence $\left(r(k)\right)_{k \ge 0}$
appears as Number A006893 in \cite{S--EIS},
where it is described as counting
trees of a kind which can be put in bijective correspondence
with the conjugacy classes of $W(C_2,k)$:
 $$
r(0) = 1, \hskip.5cm
r(1) = 2, \hskip.5cm
r(2) = 5, \hskip.5cm
r(3) = 20, \hskip.5cm
r(4) = 230, \hskip.5cm
r(5) = 26 \, 795, \hskip.5cm \hdots
$$

\medskip

To conclude this section, let us record a finite group analogue
of the formula of Witten quoted in our introduction:
if $S_n$ denotes the symmetric group on $n$ letters,
$g$ an integer, $g \ge 2$,
and $\Gamma_g$ the fundamental group 
of an oriented closed surface of genus $g$,
then the number of homomorphisms of $\Gamma_g$ to $S_n$ is given by
$$
\left\vert \operatorname{Hom} (\Gamma_g, S_n) \right\vert
\, = \, 
(n!)^{2g-2} \zeta(S_n, 2g-2) .
$$
For discussion and references, see Chapter~14 in \cite{LuSe--03}.

\head{\bf
10.~A unitary variation 
\\
and a question 
on some strengthening of Property~(T)
}\endhead

Let $G$ be a topological group. 
For all $n \ge 1$, denote by $u_n(G)$ the number (up to equivalence) 
of irreducible {\it unitary} representations of $G$
in the Hermitian space $\bold C^n$.
If $u_n(G) < \infty$ for all $n \ge 1$, set
$$
\zeta^{(u)}_G(s) \, = \, \sum_{n \ge 1} u_n(G) n^{-s} .
$$
Of course, $\zeta^{(u)}_G(s) = \zeta_G(s)$ if $G$ is compact,
and in many other cases, but not always; 
indeed, there are groups which are not rigid but are such that
$u_n(G) < \infty$ for all $n \ge 1$; we show below that
$$
\Gamma \, = \,  SL_3(\bold Z [X])
$$
is such an example.

Define for any complex number $z \in \bold C$
the $3$--dimensional representation 
$$ \pi_z \, : \, \Gamma \, \longrightarrow GL_3(\bold C),
\hskip.5cm
g \, \longmapsto g(z)
$$
where $g(z)$ denotes the result of evaluating $X$ at $z$.
As the restriction of $\pi_z$ to $SL_3(\bold Z)$
is the tautological $3$--dimensional representation of $SL_3(\bold Z)$,
the representation $\pi_z$ is clearly irreducible. 
By computing the character at a well--chosen element,
for example by computing
$$
\operatorname{trace} \left( \pi_z \left( \matrix
1+X & X & 0 \\ 1 & 1 & 0 \\ 0 & 0 & 1 \endmatrix \right)\right)
=
\operatorname{trace} \left( \matrix
1+z & z & 0 \\ 1 & 1 & 0 \\ 0 & 0 & 1 \endmatrix \right)
= 3+z ,
$$
we see that the uncountably many representations $\pi_z$
are pairwise non--equivalent.
Hence $r_3(\Gamma) = \infty$, and in particular $\Gamma$
is not rigid.

Yet $\Gamma$ has Kazhdan's Property (T),
by a recent result of Leonid Vaserstein \cite{Vase};
see also \cite{Shal--06}.
It is known that $u_n(\Gamma) < \infty$ for
any finitely generated group $\Gamma$ which has Kazhdan's Property (T);
see Proposition IV of \cite{HaRV--93},
which builds up on Propositions 2.5 and 2.6 of \cite{Wang--75},
equivalently on Corollary~2 of \cite{Wass--91}.
In other words, Property (T) implies \lq\lq unitary rigidity\rq\rq .

Summing up, $\zeta^{(u)}_{\Gamma}(s)$ is well defined,
and $\zeta_{\Gamma}(s)$ is not.
There are several strengthenings of Kazhdan's Property (T),
some already existing and probably some more to come.
For countable groups,
will one of them imply {\it bona fide} rigidity?

\enddocument